\documentclass[11pt]{amsart}
\usepackage{amsmath,amssymb,latexsym,dsfont}
\usepackage[left=2.5cm,right=2.5cm,top=2.2cm,bottom=2.2cm]{geometry}

\usepackage{amsmath,amsfonts,amssymb,epsfig, textcomp}
\usepackage{graphicx,tikz}
\usepackage{hyperref, cleveref} 
\usepackage{cases,listings}
\usepackage{color}
\usepackage{mathabx}

\newtheorem{theorem}{Theorem}[section]
\newtheorem{lemma}{Lemma}[section]

\newtheorem{proposition}{Proposition}[section]
\newtheorem{corollary}{Corollary}[section]
\newtheorem{remark}{Remark}[section]

\newtheorem{definition}{Definition}[section]
\numberwithin{equation}{section}

      \newcommand{\hy}{\hat y}

   \newcommand{\tg}{\widetilde g}
   \newcommand{\ty}{\widetilde y}

      \newcommand{\N}{\mathbb{N}}

      \newcommand{\eps}{\varepsilon}
      \newcommand{\mR}{\mathbb{R}}

      \newcommand{\dsp}{\displaystyle}

      \makeatletter
      \def\@setcopyright{}
      \def\serieslogo@{}
      \makeatother
   \newcommand{\tr}{^\mathsf{T}}

\newcommand{\cL}{\mathcal L}

\newcommand{\tvarphi}{\widetilde \varphi}

\newcommand{\tlambda}{\widetilde \lambda}

\newcommand{\tz}{\widetilde z}

\newcommand{\be}{\begin{equation}}
\newcommand{\ee}{\end{equation}}

\newcommand{\cF}{{\mathcal F}}

\newcommand{\tf}{{\widetilde f}}

\newcommand{\tA}{{\widetilde A}}

\newcommand{\cD}{{\mathcal D}}

\newcommand{\mH}{\mathbb{H}}

\newcommand{\mU}{\mathbb{U}}

\newcommand{\mX}{\mathbb{X}}

\newcommand{\tR}{{\widetilde R}}

\newcommand{\tS}{{\widetilde S}}

\newcommand{\homega}{\hat \omega}

\newcommand{\tM}{\widetilde M}

\title[Stabilization of control systems]{Stabilization of control systems associated with a strongly continuous group}

\author[H.-M. Nguyen]{Hoai-Minh Nguyen}
\address[H.-M. Nguyen]{Laboratoire Jacques Louis Lions, \newline\indent
Sorbonne Universit\'e\newline\indent
Paris, France}
\email{hoai-minh.nguyen@sorbonne-universite.fr}

\begin{document}

\maketitle 
\begin{abstract} This paper is devoted to the stabilization of a linear control system $y' = A y + B u$  and its suitable non-linear variants where $(A, \cD(A))$ is an infinitesimal generator of a strongly continuous {\it group} in a Hilbert space $\mH$, and $B$ defined in a Hilbert space $\mU$ is an admissible control operator with respect to the semigroup generated by $A$. Let $\lambda  \in \mR$ and assume that, for some {\it positive} symmetric, invertible  $Q = Q(\lambda) \in \cL(\mH)$,  for some {\it non-negative}, symmetric  $R = R(\lambda) \in \cL(\mH)$, and for some {\it non-negative}, symmetric $W = W(\lambda) \in \cL(\mU)$,  it holds 
$$
A Q + Q A^* - B W B^* + Q R Q + 2 \lambda Q = 0.       
$$ 
We then present a new approach to study the stabilization of such a system and its suitable nonlinear variants. Both the stabilization using dynamic feedback controls and the stabilization using static feedback controls in a weak sense are investigated. To our knowledge, the nonlinear case is out of reach previously when $B$ is unbounded for both types of stabilization. 
\end{abstract}

\tableofcontents



\section{Introduction}

In this paper, we study the stabilization of a linear control system associated with a strongly continuous {\it group} and its related nonlinear systems. Let $\mH$ and $\mU$ be two Hilbert spaces which denote the state space and the control space, respectively. The corresponding scalar products are  $\langle \cdot, \cdot \rangle_{\mH}$ and $\langle \cdot, \cdot \rangle_{\mU}$, and the corresponding norms are $\| \cdot \|_{\mH}$ and $\| \cdot \|_{\mU}$. Let $\big(S(t) \big)_{t \in \mR} \subset \cL(\mH)$ be  a strongly continuous {\it group} on $\mH$, i.e., 
$$
S(0) = Id \, (\mbox{the identity}), 
$$
$$
S(t_1 + t_2) = S(t_1) \circ S(t_2) \quad  \forall t_1, t_2 \in \mR, 
$$
and 
$$
\lim_{t \to 0} S(t) x = x \quad \forall x \in \mH. 
$$
Here and in what follows, for two Hilbert spaces $\mX_1$ and $\mX_2$, we denote $\cL(\mX_1, \mX_2)$ the Banach space of all bounded linear applications from $\mX_1$ to $\mX_2$ with the usual norm, and we simply denote $\cL(\mX_1, \mX_1)$ by $\cL(\mX_1)$.  

Let $(A, \cD(A))$ be the infinitesimal generator of $\big(S(t) \big)_{t \in \mR}$ and denote $S(t)^*$ the adjoint of $S(t)$ for $t \in \mR$. Then $\big( S(t)^* \big)_{t \in \mR}$ is also a strongly continuous {\it group} of continuous linear operators and its infinitesimal generator is 
$(A^*, \cD(A^*))$, which is the adjoint of $(A, \cD(A))$. As usual, we equip the domain $\cD(A^*)$ with the scalar product 
$$
\langle z_1, z_2 \rangle_{\cD(A^*)} = \langle z_1, z_2 \rangle_{\mH} +  \langle A^* z_1, A^* z_2 \rangle_{\mH} \mbox{ for } z_1, z_2 \in \cD(A^*).  
$$
Then $\cD(A^*)$ is a Hilbert space. Denote $\cD(A^*)'$ the dual space of $\cD(A^*)$ with respect to $\mH$. Then 
$$
\cD(A^*) \subset \mH \subset \cD(A^*)'. 
$$ 
Let 
$$
B \in \cL(\mU, \cD(A^*)'). 
$$

In this paper, we consider the following control system, for $T>0$,  
\be \label{CS}
\left\{\begin{array}{c}
y' = Ay + Bu \mbox{ for } t \in (0, T), \\[6pt]
y(0) = y_0, 
\end{array} \right. 
\ee
where, at time $t$, the control is $u(t) \in \mU$ and the state is $y(t) \in \mH$, and $y_0 \in \mH$ is an initial datum. This control setting is standard and used to model many control systems, see, e.g., \cite{Coron07, TW09}. Interesting aspects of the controllability and the stability of \eqref{CS} can be found in \cite{Lions71,Slemrod74,CZ95,WZ98,EN00,Coron07,Zabczyk08,TW09,TWX20} and the references therein.

As usual, see, e.g.,  \cite{Coron07, TW09}, we assume that $B$ is an {\it admissible} control operator with respect to the semi-group $\big( S(t)\big)_{t \ge 0}$ in the sense that, for all $u \in L^2([0, T]; \mU)$, it holds that 
\be \label{cond-admissibility}
\varphi \in C([0, T]; \mH) \mbox{ where } \varphi(t): = \int_0^t S(t - s) B u(s) \, ds. 
\ee
As a consequence of the closed graph theorem, see e.g., \cite{Brezis-FA}, one has 
\be \label{cond-admissibility-a}
\| \varphi\|_{C([0, T]; \mH)} \le C_T \| u \|_{L^2((0, T); \mU)}. 
\ee

Let $\lambda  \in \mR$ and assume that, for some {\it positive}, symmetric, invertible  $Q = Q(\lambda) \in \cL(\mH)$, for some {\it non-negative}, symmetric $R = R(\lambda) \in \cL(\mH)$, and for some  {\it non-negative}, symmetric $W = W(\lambda) \in \cL(\mU)$, it holds 
\be \label{identity-Op}
A Q + Q A^* - B W B^* + QRQ + 2 \lambda Q = 0,   
\ee 
where \eqref{identity-Op} is understood in the following sense 
\be \label{identity-Op-meaning} 
\langle Qx, A^*y \rangle_{\mH} + \langle A^*x, Q y \rangle_{\mH}  - \langle W B^*x,  B^*y \rangle_{\mU}  + \langle R Qx, Q y \rangle_{\mH} + 2 \lambda \langle Q x, y \rangle_{\mH}= 0  \quad \forall \, x, y \in \cD(A^*). 
\ee
In this paper, given a Hilbert space $\widetilde \mH$ and an operator $\tR \in \cL (\widetilde \mH)$ being symmetric, one says that $\tR$ is non-negative if 
$$
\langle \tR x, x \rangle_{\widetilde \mH} \ge 0 \mbox{ for all } x \in \widetilde \mH, 
$$
and one says that $\tR$ is positive if, for some positive constant $C$, it holds \footnote{Thus positivity here means coercivity.} 
$$
\langle \tR x, x \rangle_{\widetilde \mH} \ge C \| x\|_{\widetilde \mH}^2  \mbox{ for all } x \in \widetilde \mH.   
$$

Recall that system \eqref{CS} is called to be exactly controllable in some positive time $T$ if for all $y_0, y_T \in \mH$, there exists $u \in L^2((0, T); \mU)$ such that 
$$
y(T) = y_T, 
$$
where $y$ is the unique weak solution of \eqref{CS} (the definition of the weak solutions is recalled in \Cref{sect-pre}). In this case, we also call that the pair $(A, B)$ is exactly controllable in some positive time $T$. It is known that \eqref{CS} is exactly controllable in time $T>0$ if and only if the following observability inequality holds, see e.g., \cite{Coron07,TW09},  
\be \label{observability-inequality-CS}
\int_0^T \|B^*e^{s A^*} x \|_{\mU}^2 \, ds \ge C \| x\|_{\mH}^2 \mbox{ for all } x \in \mH, 
\ee
where $C$ is a positive constant independent of $x$. Here and in what follows, if $\tA$ is the infinitesimal generator of the semigroup $\big(\tS(t) \big)_{t \ge 0}$ in a Hilbert space $\widetilde \mH$, we also denote $\tS(t)$ by $e^{t \tA}$ for $t \ge 0$.

Several cases of identity \eqref{identity-Op} and their associated stabilization results appeared in the linear quadratic optimal control theory \cite{FLT88} (see also \cite{PrS87,Flandoli87-LQR,LT91,WR00,Zwart96,Staffans05} and the references therein)  under assumptions that are discussed now.  Given a {\it non-negative}, symmetric $R \in \cL(\mH)$, consider the cost function 
\be
J_T (u, y) = \int_0^T \langle R y, y \rangle_{\mH}(s) + \langle u, u \rangle_{\mU} (s) \, ds  \mbox{ for } T \in (0, + \infty]. 
\ee
For $0< T < + \infty$, let $P_T \in \cL(\mH)$ be symmetric and satisfy
$$
\langle P_T y_0, y_0 \rangle_{\mH} =  \inf_{u \in L^2((0, T), \mU)} J_T (u, y),
$$ 
where $y$ is the weak solution of \eqref{CS} corresponding to $u$.  
Assume that the finite cost condition holds, i.e., 
$$
\inf_{u \in L^2((0, + \infty), \mU)} J_\infty (u, y) < + \infty, 
$$
for all $y_0 \in \mH$. Let $u_{opt}$ and $y_{opt}$ be the unique solution corresponding to the minimizing problem $\inf_{u \in L^2((0, + \infty), \mU)} J_\infty (u, y)$, i.e.,  
\be
J_\infty (u_{opt}, y_{opt}) = \inf_{u \in L^2((0, + \infty), \mU)} J_\infty (u, y), 
\ee
where $y$ is the weak solution of \eqref{CS}. Define 
\be\label{stabilisation-W1}
S_{opt} (t) y_0 = y_{opt}(t). 
\ee
Then 
\be
S_{opt} (t) y_0 = S(t) (y_0) + \int_0^t S(t-s) B u_{opt}(s) \, ds \mbox{ for } t \ge 0. 
\ee
Let $\big(A_{opt}, \cD(A_{opt})\big)$ be the infinitesimal generator of $\big(S_{opt}(t) \big)_{t \ge 0}$. 
Then the pointwise limit of $P_T$ as $T \to + \infty$ exists. Denote this limit by $P_\infty$. It follows that $P_\infty: \cD(A_{opt}) \to \cD(A^*)$ and 
\be\label{stabilisation-W2}
u_{opt}(t) = - B^* P_\infty  y_{opt}(t) \quad \mbox{ if }  \quad y_0 \in \cD(A_{opt}).  
\ee
Assume also that $R$ is invertible. Then 
\be\label{stabilisation-W3}
\mbox{$\big(S_{opt}(t) \big)_{t \ge 0}$ is exponentially stable.} 
\ee
Assertions \eqref{stabilisation-W1}-\eqref{stabilisation-W3} thus give the stabilization of \eqref{CS} by static feedback controls in a weak sense since $-B^* P_\infty$ is not defined for every element in $\mH$ when $B$ is not bounded or equivalently when $B^*$ is not bounded. 
Assume in addition that $(S(t))_{t \in \mR} $ is a group,  and $(A^*, R^{1/2})$  and  $(A, B)$ are exactly controllable in some positive time. Then  $P_\infty$ is invertible, and $Q_\infty := P_\infty^{-1}$
satisfies the dual algebraic Riccati equation 
\be \label{DARE}
A Q_\infty + Q_\infty A^* +  Q_\infty R  Q_\infty -  B B^* = 0   
\ee
in the sense 
\be \label{DARE-meaning}
\langle  Q_\infty x, A^* z \rangle_{\mH} + \langle A^* x,  Q_\infty z \rangle_{\mH} + \langle R Q_\infty x, Q_\infty z \rangle_{\mH} = \langle B^*  x, B^*  z \rangle_{\mU} \quad \mbox{ for all } x, z \in \cD(A^*).     
\ee
Identity \eqref{DARE} is a special case of \eqref{identity-Op} for which $W = I$ and $\lambda =0$. 

We have briefly mentioned so far known stabilization results related to \eqref{identity-Op} from the optimal control theory.  We next discuss quickly known results related to \eqref{identity-Op} that come from Gramian operators and are also related to the optimal control theory. Let $\tlambda > 0$ and assume that system \eqref{CS} is exactly controllable in time $T>0$.  Thus \eqref{observability-inequality-CS} holds. 
Set, with $T_* = T + \frac{1}{2 \tlambda}$,  
\be \label{choice-elambda}
e(s) = \left\{\begin{array}{cl}
 e^{- 2 \tlambda s} & \mbox{ in } [0, T], \\[6pt]
 2 \tlambda e^{- 2 \tlambda T} (T_* - s) & \mbox{ in } (T, T_*]. 
\end{array} \right. 
\ee
It is showed in \cite{Komornik97} (see also \cite{Vest13}) that \eqref{identity-Op} holds for $\lambda =0$, $W$ being the identity, and for $Q \in \cL(\mH)$ being defined by 
\be\label{Q-Kormonik}
\langle Q x_1, x_2 \rangle_{\mH} = \int_0^{T_*} e(s) \langle B^*e^{-s A^*} x_1, B^*e^{-s A^*} x_2 \rangle_{\mU} \, ds,   
\ee
and for $R \in \cL(\mH)$ being symmetric and defined by 
$$
\langle R Q x, Q x \rangle_{\mH} = - \int_0^{T_*} e'(s) \| B^* e^{-s A^*} \|_{\mU}^2 \, ds. 
$$
Previous results when $B$ is bounded were due to Slemrod \cite{Slemrod74}.  These works are inspired by the ones of Lukes \cite{Lukes68} and Kleinman \cite{Kleinman70} where the Gramian operators were introduced in the finite-dimensional setting.  In \cite{Urquiza05}, Urquiza observed in the case $A$ is skew-adjoint and $\tlambda >0$ that \eqref{identity-Op} holds for $W$ being identity, for $\lambda =0$,  $Q$ being defined by 
\be \label{choice-elambda-U}
\langle Q x_1, x_2 \rangle_{\mH} = \int_0^{\infty} e^{- 2 \tlambda s} \langle B^*e^{-s A^*} x_1, B^*e^{-s A^*} x_2 \rangle_{\mU} \, ds,   
\ee
and for $R = 2 \tlambda Q^{-1}$.  The result of Urquiza was inspired by the Bass method previously discussed by Russell \cite[page 114-115]{Russell79} following \cite[Section 10.3]{Coron07}.  In the settings of Komornik and Urquiza, one can check that 
$$
Q \mbox{ is invertible}  \mbox{ and } (A^*, R^{1/2}) \mbox{ is exactly controllable}.  
$$
One can then apply the linear quadratic optimal control theory to conclude that system \eqref{CS} is stabilizable by static feedback controls in the weak sense \eqref{stabilisation-W2}. Komornik also proved that \eqref{CS} is stabilizable with the rate $\tlambda$ and  Urquiza \cite{Urquiza05} also established that \eqref{CS} is stabilizable with the rate $2 \tlambda$ when $A$ is skew-adjoint, both are in the weak sense. To our knowledge, these known results mentioned have not been successfully extended to the nonlinear case.

\medskip
The goal of this paper is to present a new method to study the stabilization of \eqref{CS}  and its suitable nonlinear variants under condition \eqref{identity-Op}. We study the stabilization of 
\eqref{CS} by dynamic feedback controls and by static feedback controls in a weak sense, which we call a trajectory sense. 
A system is called dynamically stabilizable if it can be embedded as a subsystem of a larger, exponentially stable well-posed system.  This definition has been used for finite dimensions, see e.g.,  \cite[chapter 11]{Coron07}, and for linear systems in infinite dimension, see e.g.,  \cite{WR00}.

Our approach is essentially based on the construction of new auxiliary dynamics for both types of stabilization (see \Cref{thm-main1} and \Cref{thm-main2}) and  ``integration by parts arguments" (see \Cref{lem-transposition} and \Cref{lem-identity-Op}). The new adding variable is inspired by the adjoint state in the linear quadratic optimal control theory and the way to choose controls in the Hilbert Uniqueness Method (HUM) principle. 
The advantage of our approach is at least twofold. First, the method works well in both linear and {\it nonlinear} settings. Second, a Lyapunov function is also provided for the static feedback controls. To our knowledge, the stabilization of such systems by dynamic feedback controls is new even in the linear setting. The nonlinear case is out of reach previously when $B$ is unbounded for both types of stabilization. Concerning the static feedback controls, as far as we know, a Lyapunov function is not known even in the case where $B$ is bounded and $A$ is not;  a Lyapunov function was previously given in the finite-dimensional case \cite{Coron07, Kleinman70}.
Consequently, we derive that if the system is exactly controllable in some positive time, then the system is rapidly stabilizable. The techniques and ideas used in this paper have been applied and combined with the ideas in \cite{CoronNg17} to study the finite-time stabilization of the Schr\"odinger equation with bilinear controls \cite{Ng-S-stabilization} and of the KdV equations \cite{Ng-KdV-stabilization}. 

Adding a new variable is very natural and has been used a long time ago in the control theory even in finite dimensions for linear control systems, see e.g.,  \cite[Section 11.3]{Coron07} and \cite[Chapter 7]{Sontag98}.  Coron and Pradly \cite{CP91} showed that there exists a nonlinear system in finite dimensions for which the system cannot be stabilized by static feedback controls but can be stabilized by dynamic feedback ones.  
Dynamic feedback controls of finite dimensional nature, i.e., the complement system is a system of differential equations, have been previously implemented in the infinite dimensions, see e.g., \cite{CVKB13,Coron-Ng22}.  It is interesting to know whether or not adding a new variable is necessary in the setting of this paper.

\medskip 
The paper is organized as follows. In \Cref{sect-results}, we state the main results of the paper on the dynamic feedback and the static feedback in the trajectory sense. \Cref{sect-pre} is devoted to the well-posedness and some properties of various linear systems considered in this paper. The proofs of the main results on the dynamic feedback and the static feedback are given in  \Cref{sect-dynamic-feedback} and \Cref{sect-static-feedback}, respectively. In \Cref{sect-generator}, we also discuss the infinitesimal generator of the semigroup associated with the static feedback controls given in \Cref{thm-main1}, this in particular implies new information on $(A_{opt}, \cD(A_{opt}))$.  
Finally, in \Cref{sect-Q-E}, we discuss choices of $Q$ (and also $R$ and $W$) when the system is exactly controllable.  

\section{Statement of the main results} \label{sect-results}

This section consisting of two subsections is organized as follows. In the first  
subsection, we discuss the stabilization \eqref{CS} by dynamic feedback controls.  In the second subsection, we discuss the stabilization of \eqref{CS} by static feedback controls in the trajectory sense.  Here and in what follows in this section, we always assume that $(A, \cD(A))$ is an infinitesimal generator of a strongly continuous {\it group} in $\mH$, and $B \in \cL(\mU, \cD(A^*)')$ is an admissible control operator with respect to the semigroup generated by $A$.

\subsection{Stabilization by dynamic feedback controls}

Given an infinitesimal generator $\tA$ of a semigroup in a Hilbert space $\widetilde \mH$, set 
$$
\omega_0(\tA) = \inf_{t > 0} \log \| e^{t \tA} \|_{\cL(\widetilde \mH)},  
$$
which denotes the growth of the $e^{t \tA}$ for $t \ge 0$. It is known, see e.g., \cite{EN00}, that 
$$
-\infty \le \omega_0(\tA) < + \infty. 
$$
Concerning the dynamic feedback controls of \eqref{CS}, we have the following result. 

\begin{theorem}\label{thm-main1-D} Let $\lambda \in \mR$ and assume \eqref{identity-Op} with $R = 0$, and let $\lambda_1 \in \mR$. Let $\homega_0(A) \ge \omega_0(A)$ and $\homega_0(-A^*) \ge \omega_0(-A^*)$ be two real constants such that, for some positive constant $c$,  
$$
\|e^{t A}\|_{\cL(\mH)} \le c e^{t \homega_0(A)} \mbox{ for } t \ge 0
\quad \mbox{ and } \quad 
\|e^{-t A^*}\|_{\cL(\mH)} \le c e^{t \homega_0(-A^*)} \mbox{ for } t \ge 0,     
$$
and assume that 
\be \label{thm-main1-lambda-D}
\lambda_1 - 2 \lambda  >  \homega_0 (A) - \homega_0(-A^*).       
\ee
Given $y_0, \ty_0 \in \mH$ arbitrary, let $(y, \ty)\tr \in \big( C^0([0, T]; \mH) \big)^2$ be the unique weak solution of the system 
\be \label{thm-main1-sys}
\left\{\begin{array}{c}
y' = A y - B W B^* \ty \quad \mbox{ in } (0, + \infty), \\[6pt]
\ty'  = - A^* \ty - 2 \lambda \ty + \lambda_1 Q^{-1} (y - Q\ty) \quad \mbox{ in } (0, + \infty), \\[6pt]
y(0) = y_0, \quad  \ty(0) = \ty_0. 
\end{array} 
\right. 
\ee
Then 
\be \label{thm-main1-cl1-D}
\| y(t)\|_{\mH} + \| \ty(t)\|_{\mH} \le C  e^{(\homega_0(-A^*) - 2 \lambda) t}  \big( \| y(0) \|_{\mH} + \|\ty(0)\|_{\mH} \big) \mbox{ for } t \ge 0, 
\ee
where $C$ is a positive constant independent of $t$ and $(y_0, \ty_0)$.  Consequently, if $A$ is skew-adjoint and $\lambda_1 > 2 \lambda$,  then  
\be \label{thm-main1-cl2-D}
\| y(t)\|_{\mH} + \| \ty(t)\|_{\mH} \le C  e^{- 2 \lambda t}  \big( \| y(0) \|_{\mH} + \|\ty(0)\|_{\mH} \big) \mbox{ for } t \ge 0.  
\ee
\end{theorem}

\begin{remark} \rm
The well-posedness of the weak solutions in \Cref{thm-main1-D} is established in \Cref{lem-WP-thm-main1-D}. 
\end{remark}

We next illustrate how this result can be extended to a nonlinear setting. Let $f: \mH \to \mH$ be continuous such that for all $\eps > 0$ there exists $\delta > 0$ such that 
\be \label{cond-NL1}
\|f(x)\|_{\mH} \le \eps \|x \|_{\mH} \mbox{ for } x \in \mH \mbox{ with } \|x\|_{\mH} < \delta,  
\ee
and $f$ is Lipschitz in a neighborhood of $0$ in $\mH$, i.e., there exist $r > 0$ and $\Lambda >0$ such that 
\be \label{cond-NL2}
\|f(x) - f(y)\|_{\mH} \le \Lambda \|x -y\|_{\mH} \mbox{ for } x, \, y \in \mH \mbox{ with } \| x\|_{\mH}, \, \| y\|_{\mH} < r. 
\ee
We consider the following control system 
\be \label{CS-NL}
\left\{\begin{array}{c}
y' = Ay + f (y) + Bu \mbox{ for } t \in (0, T), \\[6pt]
y(0) = y_0 \in \mH. 
\end{array} \right. 
\ee

Concerning the local stabilization of \eqref{CS-NL}, we have the following stabilization result. 

\begin{theorem}\label{thm-main2-D} Let $\lambda \in \mR$  and assume that \eqref{identity-Op} holds with $R =0$, and let $\lambda_1, \gamma  \in \mR$ be such that $\gamma < \lambda$. Let $\homega_0(A) \ge \omega_0(A)$ and $\homega_0(-A^*) \ge \omega_0(-A^*)$ be two real constants such that, for some positive constant $c$,  
$$
\|e^{t A}\|_{\cL(\mH)} \le c e^{t \homega_0(A)} \mbox{ for } t \ge 0
\quad \mbox{ and } \quad 
\|e^{-t A^*}\|_{\cL(\mH)} \le c e^{t \homega_0(-A^*)} \mbox{ for } t \ge 0.      
$$
Assume that 
\be \label{thm-main2-lambda-D}
\lambda_1 - 2 \lambda  >  \homega_0 (A) - \homega_0(-A^*),  \quad 2 \gamma - \homega_0(-A^*) > 0,        
\ee
and \eqref{cond-NL1} and \eqref{cond-NL2} hold. There exists $\eps > 0$ (small) such that for $y_0, \ty_0 \in \mH$ with $\| (y_0, \ty_0)\tr\|_{\mH}  \le \eps$, there exists a unique solution $(y, \ty)\tr \in \big( C^0([0, T]; \mH) \big)^2$ of the system 
\be \label{thm-main2-sys-D}
\left\{\begin{array}{c}
y' = A y + f (y) - B W B^* \ty \quad \mbox{ in } (0, + \infty), \\[6pt]
\ty'  = - A^* \ty - 2 \lambda \ty  + Q^{-1} f (Q \ty) +  \lambda_1 Q^{-1} (y - Q\ty) \quad \mbox{ in } (0, + \infty), \\[6pt]
y(0) = y_0, \quad  \ty(0) = \ty_0. 
\end{array} 
\right. 
\ee
Moreover, we have  
\be \label{thm-main2-cl1-D}
\| y(t)\|_{\mH} + \| \ty(t)\|_{\mH} \le C  e^{(\homega_0(-A^*) - 2 \gamma) t}  \big( \| y(0) \|_{\mH} + \|\ty(0)\|_{\mH} \big) \mbox{ for } t \ge 0, 
\ee
where $C$ is a positive constant independent of $t$ and $(y_0, \ty_0)$. 
Consequently, if $A$ is skew-adjoint and $\lambda_1 > 2 \lambda  > 2 \gamma > 0$
then  
\be \label{thm-main2-cl3}
\| y(t)\|_{\mH} + \| \ty(t)\|_{\mH} \le C  e^{- 2 \gamma t}  \big( \| y(0) \|_{\mH} + \|\ty(0)\|_{\mH} \big) \mbox{ for } t \ge 0. 
\ee
\end{theorem}

\begin{remark} \rm The weak solutions given in \Cref{thm-main2-D} are understood in the sense of the weak solutions where the nonlinear terms play 
as a part of the source term. 
\end{remark}

\begin{remark} \rm
The well-posedness of the weak solutions in \Cref{thm-main2-D} is a part of the proof. In comparison with \Cref{thm-main1-D}, $\lambda$ is supposed to satisfy the condition $2 \lambda - \homega_0(-A^*) > 0$ in \Cref{thm-main2-D} to make sure that the solution remains small for large time. 
\end{remark}

As a consequence of \Cref{thm-main1-D} and \Cref{thm-main2-D} (see also \Cref{pro-Q-E-U}), we obtain the following results. 

\begin{proposition} Assume that system \eqref{CS} is exactly controllable in some positive time. System \eqref{CS} is rapidly dynamically stabilizable. 
\end{proposition}

\begin{proposition} Assume that system \eqref{CS} is exactly controllable in some positive time, and  \eqref{cond-NL1} and \eqref{cond-NL2} hold. 
System \eqref{CS-NL} is locally rapidly dynamically stabilizable. 
\end{proposition}

Recall that system \eqref{CS} is called rapidly dynamically stabilizable if it can be dynamically exponentially stabilizable with an arbitrary decay rate. A similar meaning with suitable modifications is used for system \eqref{CS-NL}.  

\subsection{Stabilization by static feedback controls} Here is the first main result on the static feedback controls of \eqref{CS}. 

\begin{theorem}\label{thm-main1} Let $\lambda \in \mR$ and assume \eqref{identity-Op}.  Given $y_0 \in \mH$, let $(y, \ty)\tr \in \big( C^0([0, T]; \mH) \big)^2$ be the unique weak solution of the system 
\be \label{thm-main1-sys}
\left\{\begin{array}{c}
y' = A y - B W B^* \ty \quad \mbox{ in } (0, + \infty), \\[6pt]
\ty'  = - A^* \ty - 2 \lambda \ty - R Q \ty \quad \mbox{ in } (0, + \infty), \\[6pt]
y(0) = y_0, \quad  \ty(0) = \ty_0: = Q^{-1} y_0. 
\end{array} 
\right. 
\ee
Then 
\be \label{thm-main1-cl1}
\ty(t) = Q^{-1}y (t)  \mbox{ for } t \ge 0, 
\ee
and  
\begin{multline} \label{thm-main1-cl2}
\| Q^{-1/2} y (t)\|_{\mH}^2 - \| Q^{-1/2} y (\tau)\|_{\mH}^2
 \\[6pt]
= - 2 \lambda \int_\tau^t \|Q^{-1/2} y(s) \|_{\mH}^2 \, ds  - \int_\tau^t \Big( \|W^{1/2} B^* \ty (s) \|_{\mU}^2 + \| R^{1/2} y (s)\|_{\mH}^2 \Big) \, d s \mbox{ for } t \ge \tau \ge 0.   
\end{multline}
Consequently, 
\be \label{thm-main1-cl3}
\| Q^{-1/2} y (t)\|_{\mH}  \le e^{- \lambda t} \| Q^{-1/2} y (0)\|_{\mH}   
\mbox{ for } t \ge 0.   
\ee
\end{theorem}

Some comments on \Cref{thm-main1} are in order. Since 
$$
\ty'  = - A^* \ty - 2 \lambda \ty - R Q \ty \mbox{ in } (0, + \infty)
$$
and $\ty(0) \in \mH$, it follows from \Cref{lem-transposition} given in \Cref{sect-pre}  that $\ty \in C([0, T); \mH)$ is well-defined for all $T > 0$ and moreover, 
$$
B^*\ty \in L^2((0, T), \mH) \mbox{ for all } T > 0. 
$$
We thus derive that system~\eqref{thm-main1-sys} is well-posed and \eqref{thm-main1-cl2} makes sense. Combing \eqref{thm-main1-cl1} and the equation of $y$
$$
y' = A y - B W B^* \ty, 
$$
we have thus shown that the control system $y' = A y + B u$ with the {\it static feedback} control 
\be
\label{thm-main1-uuu}
``u = - W B^* Q^{-1} y" \mbox{ for } t \ge 0, 
\ee is well-posed in the sense given in \Cref{thm-main1}.  We only consider \eqref{thm-main1-uuu} as static feedback controls in a weak sense, which we call {\it a trajectory sense},  since for $y \in \mH$, it is not clear how to give the sense to the action  $-W B^* Q^{-1} y$. In comparison with the static feedback controls in the sense given by \eqref{stabilisation-W2}, the static feedback controls given \eqref{thm-main1-uuu} are well-defined in the sense of \Cref{thm-main1} for all initial data $y_0 \in \mH$. \Cref{thm-main1} can be considered as a new way to view the feedback controls given in \eqref{stabilisation-W2}.

It is important in \Cref{thm-main1} that $\ty_0 = Q^{-1} y_0$ in \eqref{thm-main1-sys}. Due to this fact, one cannot derive from \Cref{thm-main1} that system \eqref{CS} is dynamically stabilizable via the system 
\be
\left\{\begin{array}{c}
y' = A y - B W B^* \ty \quad \mbox{ in } (0, + \infty), \\[6pt]
\ty'  = - A^* \ty - 2 \lambda \ty - R Q \ty \quad \mbox{ in } (0, + \infty). 
\end{array} 
\right. 
\ee
This is the reason to introduce the term $\lambda_1 Q^{-1} (y - Q \ty)$ in \Cref{thm-main1-D}. 

\begin{remark}\rm From \eqref{thm-main1-cl2}, the quantity $\| Q^{-1/2} y (t)\|_{\mH}^2$ can be viewed as the Lyapunov function of the system. This fact seems new to us even in the case where $B$ is bounded and $A$ is not. 
\end{remark}

\begin{remark} \rm
Assertion~\eqref{thm-main1-cl1} was known in the case where $\lambda = 0$, $W = I$, and under the additional assumptions that $(A^*, R^{1/2})$ and $(A, B)$ are exactly controllable in some positive time, see \cite[Theorems 2.4, 2.6, and 2.7]{FLT88}.  
\end{remark}

We next present a consequence of \Cref{thm-main1} in the case where $A$ is a skew-adjoint operator and $R = 0$. 

\begin{corollary}\label{cor-main1} Let $\lambda \in \mR$,  and assume that \eqref{identity-Op} holds with $R = 0$ and $A$ is skew-adjoint.  Given $y_0 \in \mH$, let $(y, \ty)\tr \in \big( C^0([0, T]; \mH) \big)^2$ be the unique weak solution of the system \eqref{thm-main1-sys} with $R = 0$. Then \eqref{thm-main1-cl1} holds and, for some positive constants $C_1, C_2$, independent of $y_0$,   
\be \label{cor-main1-cl3}
C_1 e^{-2 \lambda t} \le \| y(t) \|_{\mH}  \le C_2 e^{- 2 \lambda t} \| y_0\|_{\mH} 
\mbox{ for } t \ge 0. 
\ee
\end{corollary}

\Cref{cor-main1} is a direct consequence of \Cref{thm-main1}. Indeed, \eqref{thm-main1-cl1} is a consequence of \Cref{thm-main1}. Since $A$ is skew-adjoint, it follows from the equation of $\ty$ that 
\be \label{cor-main1-p1}
\| \ty (t) \| = e^{-2 \lambda t} \| \ty (0) \| \mbox{ for } t \ge 0. 
\ee
Assertion \eqref{cor-main1-cl3}  is now a consequence of \eqref{thm-main1-cl1} and \eqref{cor-main1-p1}.

\medskip 
We next deal with the local stabilization of \eqref{CS-NL} by static feedback controls in the trajectory sense. 

\begin{theorem}\label{thm-main2} Let $\lambda > 0$ and assume \eqref{identity-Op}, \eqref{cond-NL1}, and \eqref{cond-NL2}.  There exists $\eps > 0$ (small) such that for $y_0 \in \mH$ with $\| y_0\|_{\mH} \le \eps$, there exists a unique weak solution $(y, \ty)\tr \in \big( C^0([0, T]; \mH) \big)^2$ of the system 
\be \label{thm-main2-sys-00}
\left\{\begin{array}{c}
y' = A y + f (y) - B W B^* \ty \quad \mbox{ in } (0, + \infty), \\[6pt]
\ty'  = - A^* \ty - 2 \lambda \ty - R Q \ty + Q^{-1} f (Q \ty) \quad \mbox{ in } (0, + \infty), \\[6pt]
y(0) = y_0, \quad  \ty(0) = \ty_0: = Q^{-1} y_0. 
\end{array} 
\right. 
\ee
Moreover, we have  
\be \label{thm-main2-cl1}
\ty = Q^{-1}y \mbox{ for } t \ge 0, 
\ee
and  
\begin{multline} \label{thm-main2-cl2}
\| Q^{-1/2} y (t)\|_{\mH}^2  - \| Q^{-1/2} y (\tau)\|_{\mH}^2
= - 2 \lambda \int_\tau^t \| Q^{-1/2} y (s)\|_{\mH}^2 \, ds   \\[6pt]
- \int_\tau^t \Big(\|W^{1/2} B^* \ty (s) \|_{\mU}^2 + \| R^{1/2} y (s)\|_{\mH}^2 \Big) \, d s + 2 \int_\tau^t \langle   f(y(s)), Q^{-1} y(s) \rangle \, ds
\mbox{ for } t \ge \tau \ge 0.   
\end{multline}
Consequently, for all $0< \gamma < \lambda$, there exists $\eps_\gamma$ such that for $y_0 \in \mH$ with $\| y_0\|_{\mH} \le \eps_\gamma$, it holds 
\be \label{thm-main2-cl3}
\| Q^{-1/2} y (t)\|_{\mH}  \le e^{-  \gamma t} \| Q^{-1/2} y (0)\|_{\mH}
\mbox{ for } t \ge 0.   
\ee\end{theorem}

\begin{remark} \rm The weak solutions given in \Cref{thm-main2} are understood in the sense of the weak solutions given in \Cref{sect-pre} where the nonlinear term plays as a part of the source term. 
\end{remark}

\begin{remark} \rm In comparison with \Cref{thm-main1-D}, $\lambda$ is supposed to be positive in \Cref{thm-main2-D} to make sure that the solution remains small for large time. 
\end{remark}

Here is a variant of \Cref{cor-main2} in the nonlinear setting, which is a direct consequence of \Cref{thm-main2}, and the proof is omitted. 

\begin{corollary}\label{cor-main2} Let $\lambda > 0$,  and assume that \eqref{identity-Op} holds with $R = 0$ and $A$ is skew-adjoint.  
Assume \eqref{cond-NL1} and \eqref{cond-NL2}. There exists $\eps > 0$ (small) such that for $y_0 \in \mH$ with $\| y_0\|_{\mH} \le \eps$, there exists a unique solution $(y, \ty)\tr \in \big( C^0([0, T]; \mH) \big)^2$ of the system \eqref{thm-main2-sys-00} with $R = 0$. 
Moreover, \eqref{thm-main2-cl1} holds, and, for all $0< \gamma < \lambda$, there exists $\eps_\gamma$ such that for $y_0 \in \mH$ with $\| y_0\|_{\mH} \le \eps_\gamma$, it holds, for some positive constants $C$, independent of $y_0$,   
\be \label{cor-main2-cl3}
\| y(t) \|_{\mH}  \le C e^{- 2 \gamma t} \| y_0\|_{\mH} 
\mbox{ for } t \ge 0. 
\ee
\end{corollary}

As a consequence of \Cref{thm-main1} and \Cref{thm-main2} (see also \Cref{pro-Q-E}), we obtain the following results. 

\begin{proposition} Assume that system \eqref{CS} is exactly controllable in some positive time. System \eqref{CS} is rapidly (statically) stabilizable in the trajectory sense. 
\end{proposition}

\begin{proposition} Assume that system \eqref{CS} is exactly controllable in some positive time, and \eqref{cond-NL1} and \eqref{cond-NL2} hold.  System \eqref{CS-NL} is locally rapidly (statically) stabilizable in the trajectory sense. 
\end{proposition}

\section{Preliminaries} \label{sect-pre}

In this section, we state and prove the well-posedness and some properties of various linear control systems considered in this paper.  It is more convenient to consider a slightly more general system 
\be \label{CS-G}
\left\{\begin{array}{c}
y' = Ay + f + Bu + My  \mbox{ in } t \in (0, T), \\[6pt]
y(0) = y_0, 
\end{array} \right. 
\ee
with $y_0 \in \mH$, and $f \in L^1((0, T); \mH)$, and $M \in \cL(\mH)$. Recall that $B$ is assumed to be an admissible control operator with respect to the semigroup $\big( S(t) \big)_{t \ge 0} \subset \cL(\mH)$ generated by the operator $A$ {\it throughout} the paper.   In this section, we only assume that $\big( S(t) \big)_{t \ge 0} \subset \cL(\mH)$ is a strongly continuous {\it semigroup}.  
A weak solution $y$ of \eqref{CS-G} is understood as an element $y \in C([0, T]; \mH)$ such that 
\be\label{meaning-CS-G}
\left\{\begin{array}{c} \frac{d}{dt} \langle y, \varphi \rangle_{\mH}  = \langle A y + f + Bu + My, \varphi \rangle_{\mH} \mbox{ in } (0, T) \\[6pt]
y(0) =  y_0
\end{array}\right. \mbox{ for all } \varphi \in \cD(A^*)
\ee
for which 
\begin{itemize}
\item[$i)$] the differential equation in \eqref{meaning-CS-G} is understood in the distributional sense, 
\item[$ii)$] the term $ \langle A y + f + Bu + My, \varphi \rangle_{\mH} $ is understood as $\langle y, A^*\varphi \rangle_{\mH}  +  \langle f + My, \varphi \rangle_{\mH} + \langle u, B^*\varphi \rangle_{\mU}$. 
\end{itemize}
The convention in $ii)$ will be used throughout this section. 

\medskip 

We begin by recalling the well-posedness of \eqref{CS-G}, see \cite[Sections 4.1 and 4.2]{TW09} (in particular, \cite[Remark 4.1.2 and Proposition 4.2.5]{TW09}) \footnote{There is no $f$ in the statement of \cite[Proposition 4.2.5]{TW09} but the result also holds with $f \in L^1((0, T); \mH)$ and the analysis is the same.}.

\begin{proposition}\label{pro-WP} Let $T>0$, $y_0 \in \mH$,  $f \in L^1((0, T); \mH)$, and $M \in \cL(\mH)$.  Then  
\begin{itemize}
\item[$i)$] $y \in C([0, T], \mH)$ is a weak solution of \eqref{CS-G}
if and only if, with $\tf := f + Bu + My$, it holds \footnote{This identity is understood in $\cD(A^*)'$, i.e., $\langle y(t), \varphi \rangle_{\mH}  = \langle S(t) y_0, \varphi \rangle_{\mH} + \int_0^t \langle S(t - s) \tf(s), \varphi \rangle_{\mH} \, ds$ in  $[0, T]$ for all $\varphi \in \cD(A^*)$.}
\be\label{pro-WP-identity}
y(t) = S(t) y_0 + \int_0^t S(t - s) \tf(s) \, ds \mbox{ for } t \in [0, T]. 
\ee
\item [$ii)$] there exists a unique weak solution $y \in C([0, T], \mH)$ of \eqref{CS-G}. 
\end{itemize}
\end{proposition}

\begin{remark} \rm  The equivalence between weak solutions and mild solutions was first proved in the case $B$ is bounded and $f \in C([0, T]; \mH)$ by Ball \cite{Ball77}, see also \cite[Chapter 1 of Part II]{BDDM07} for related results when $B$ is bounded. 
\end{remark}

The unique weak solution given in \Cref{pro-WP} also satisfies the transposition meaning as established in the following result, which is one of the key technical result of this paper. 

\begin{lemma}\label{lem-transposition}
Let $T>0$, $y_0 \in \mH$,  $f \in L^1((0, T); \mH)$, and $M \in \cL(\mH)$,  and let $y \in C([0, T]; \mH)$ be the unique weak solution of \eqref{CS-G}.  We have, for $t \in (0, T]$,  for $z_t \in \cD(A^*)$, and for $g \in C([0, t]; \cD(A^*))$,  
\begin{multline}\label{lem-transposition-cl1}
 \langle y(t), z_t \rangle_{\mH}  - \langle y_0, z(0) \rangle_{\mH} =  \int_0^t \langle u (s), B^*z(s) \rangle_{\mU} \, ds  \\[6pt] -  \int_0^t \langle g(s),  y(s) \rangle_{\mH} \, ds 
+  \int_0^t \langle f(s),  z(s) \rangle_{\mH} \, ds + \int_0^t \langle M y(s), z (s) \rangle_{\mH} \, ds, 
\end{multline}
where $z \in C([0, t]; \mH)$ is the unique weak solution of the backward system 
\be \label{lem-transposition-z}
\left\{\begin{array}{c}
z' = - A^* z - g  \mbox{ in }  (0, t), \\[6pt]
z(t)  = z_t.  
\end{array} \right. 
\ee
Consequently, for $z_T \in \mH$ and $g \in L^1((0, T); \mH)$, the unique weak solution $z \in C([0, T]; \mH)$ of \eqref{lem-transposition-z} with $t = T$ satisfies
\be \label{lem-transposition-cl2}
\| B^*z \|_{L^2((0, T); \mU)} \le C_T \Big(\| g\|_{L^1((0, T); \mH)} + \| z_T\|_{\mH} \Big),
\ee
and \eqref{lem-transposition-cl1} holds for $z_t \in \mH$ and $g \in L^1((0, t); \mH)$. Here $C_T$ denotes a position constant independent of $g$, $f$,  and  $z_T$. 
\end{lemma}

\begin{remark} \label{rem-lem-transposition} \rm For $0< T \le T_0$, the constant $C_T$ in \eqref{lem-transposition-cl2} can be chosen independent of $T$. In fact, extend $g$ by $0$ for $t < 0$ and denote this extension by $\tg$. Consider the weak solution $\tz$ of the system 
\be \label{rem-lem-transposition-z}
\left\{\begin{array}{c}
\tz' = - A^* \tz - g  \mbox{ in }  (T-T_0, T), \\[6pt]
\tz(T)  = z_T.  
\end{array} \right. 
\ee
By \eqref{lem-transposition-cl2}, we have 
$$
\| B^* \tz\|_{L^2((T-T_0, T); \mU)} \le C_{T_0} \Big(\| \tg\|_{L^1((T-T_0, T); \mH)} + \| z_T\|_{\mH} \Big). 
$$
The desired assertion follows by noting that $\tz = z$ in $(0, T)$ and using the definition of $g$. 
\end{remark}

In what follows, for notational ease, we use $\langle \cdot, \cdot, \rangle$ to denote $\langle \cdot, \cdot, \rangle_{\mH}$ or $\langle \cdot, \cdot, \rangle_{\mU}$ in a clear context.   We now give the proof of \Cref{lem-transposition}. 

\begin{proof}[Proof of \Cref{lem-transposition}] Let $z_t \in \cD(A^*)$ and $g \in C([0, t]; \cD(A^*))$, and let $z \in C([0, t]; \mH)$ be the unique weak solution of \eqref{lem-transposition-z}. We have, for $n \ge 2$,  
$$
 \langle y(t), z(t) \rangle - \langle y(0), z(0) \rangle = \sum_{i=1}^n \Big(  \langle y(t_i), z(t_i) \rangle - \langle y(t_{i-1}), z(t_{i-1}) \rangle \Big),  
$$
where $t_0 = 0$ and $t_{i} = t_{i-1} + t/ n$ for $1 \le i \le n$. 

Since $z_t \in \cD(A^*)$ and $g \in C([0, t]; \cD(A^*))$, it follows that $z \in C\big([0, t]; \cD(A^*)\big)$. 
We thus obtain   
\begin{multline}\label{lem-transposition-p0-1}
\langle y(t_i), z(t_i) \rangle - \langle y(t_{i-1}), z(t_{i-1}) \rangle  = \langle y(t_i), z(t_i)-z(t_{i-1}) \rangle +  \langle y(t_i) - y(t_{i-1}), z(t_{i-1})  \rangle  \\[6pt]
\mathop{=}^{\eqref{CS-G}, \eqref{lem-transposition-z}} \langle y(t_i), \int_{t_{i-1}}^{t_i} \big(-A^* z(s) - g(s) \big) \, ds \rangle +   \int_{t_{i-1}}^{t_i} \langle A y(s) + \tf(s), z(t_{i-1}) \rangle \, ds. 
\end{multline}
where $\tf = f + Bu + My$.  Recall that the convention $ii)$ in the definition of the weak solutions of \eqref{meaning-CS-G} is used here.  Using the fact $z \in C\big([0, t]; \cD(A^*)\big)$ and $y \in C\big([0, t]; \mH \big)$, we derive that 
\begin{multline}\label{lem-transposition-p0-2}
\langle y(t_i), \int_{t_{i-1}}^{t_i} \big(-A^* z(s) - g(s) \big) \, ds \rangle +  \int_{t_{i-1}}^{t_i} \langle A y(s) + \tf(s) , z(t_{i-1}) \rangle  \, ds \\[6pt]
= \int_{t_{i-1}}^{t_i} \langle y(s),  \big(-A^* z(s) - g(s) \big) \, ds \rangle + 
 \int_{t_{i-1}}^{t_i} \langle  A y(s) + \tf(s), z(s) \rangle \, ds + o(t_i - t_{i-1}). 
\end{multline}
Here the standard notation of $o(\cdot)$ is used: $o(s)/|s| \to 0$ as $s \to 0$. Combining \eqref{lem-transposition-p0-1} and \eqref{lem-transposition-p0-2} yields 
$$
\langle y(t_i), z(t_i) \rangle - \langle y(t_{i-1}), z(t_{i-1}) \rangle  =    - \int_{t_{i-1}}^{t_i} \langle y(s),  g(s) \, ds \rangle + 
\int_{t_{i-1}}^{t_i} \langle   \tf(s) , z(s) \rangle \, ds + o(t_i - t_{i-1}). 
$$
Using the definition of $\tf$, we derive that  
\begin{multline*}
\langle y(t_i), z(t_i) \rangle - \langle y(t_{i-1}), z(t_{i-1}) \rangle  
=  \int_{t_{i-1}}^{t_i} \langle  u (s), B^* z(s) \rangle \, ds  -  \int_{t_{i-1}}^{t_i} \langle g(s),  y(s) \rangle \, ds  \\[6pt]
 +  \int_{t_{i-1}}^{t_i}\langle f(s),  z(s) \rangle \, ds + \int_{t_{i-1}}^{t_i} \langle M y(s), z (s) \rangle \, ds + o(t_i - t_{i-1}). 
\end{multline*}
Summing with respect to $n$ and letting $n \to + \infty$, we reach \eqref{lem-transposition-cl1} for $z_t \in \cD(A^*)$ and $g \in C([0, t]; \cD(A^*))$.  

\medskip 
We next deal with \eqref{lem-transposition-cl2}. Fix $z_T \in \cD(A^*)$ and $g \in C([0, T]; \cD(A^*))$. Let  $u \in L^2((0, T); \mU)$ and let $y \in C([0, T]; \mH)$ be the unique weak solution of  \eqref{CS-G} with $f=0$, $y_0 = 0$, and $M=0$.  
Applying \eqref{lem-transposition-cl1} with $t = T$, we have 
\be \label{lem-transposition-p21}
\int_0^T \langle u (s), B^*z(s) \rangle \, ds = \langle y(T), z_T \rangle +   \int_0^T \langle g(s),  y(s) \rangle \, ds. 
\ee
Since 
\begin{multline} \label{lem-transposition-p22}
|\langle y(T), z_T \rangle| +   \int_0^T |\langle g(s),  y(s) \rangle| \, ds \le \| y(T)\| \|z_T \| + \| g\|_{L^1((0, T);  \mH)} \| y\|_{L^\infty((0, T);  \mH)} \\[6pt] 
\mathop{\le}^{\eqref{cond-admissibility-a}, \Cref{pro-WP}} C \| u\|_{L^2((0, T); \mU)}  \Big( \|z_T \| + \| g\|_{L^1((0, T);  \mH)} \Big).  
\end{multline}
Combining \eqref{lem-transposition-p21} and \eqref{lem-transposition-p22} yields 
$$
\| B^*z \|_{L^2((0, T), \mU)} \le C \Big( \|z_T \| + \| g\|_{L^1((0, T);  \mH)} \Big).
$$ 
Assertion of \eqref{lem-transposition-cl2} in the case $z_T \in \mH$ and $g \in L^1((0, T); \mH)$ follows from this case by density. 

Finally, \eqref{lem-transposition-cl1} with  $z_t \in \mH$ and $g \in L^1((0, t); H)$ also follows from the case $z_t \in \cD(A^*)$ and $g \in C([0, t]; \cD(A^*))$ by density. 
\end{proof}

We now prove that the solutions in the transposition sense are also unique. Their existence is a direct consequence of \Cref{pro-WP} and \Cref{lem-transposition}. We first state the meaning of transposition solutions of system \eqref{CS-G}. 

\begin{definition}\label{def-transposition} Let $T>0$, $y_0 \in \mH$,  $f \in L^1((0, T); \mH)$, and $M \in \cL(\mH)$. A function $y \in C([0, T]; \mH)$ is called a transposition solution of \eqref{CS-G} if for all $t \in (0, T]$,  $z_t \in \mH$, and $g \in L^1((0, t); H)$, identity \eqref{lem-transposition-cl1} holds where $z \in C([0, t]; \mH)$ is the unique weak solution of \eqref{lem-transposition-z}. 
\end{definition}

We have the following result. 
 
\begin{lemma}\label{lem-WP-Trans} Let $T>0$, $y_0 \in \mH$,  and $u \in L^2((0, T); \mU)$. There exists a unique transposition solution $y \in C^0([0, T]; \mH)$ of \eqref{CS-G}. Moreover, 
\be \label{lem-WP-Trans-cl}
\| y(\tau) \|_{\mH} \le C_T \Big( \| y_0 \| + \| u\|_{L^2((0, T); \mU)} \Big), 
\ee
for some positive constant $C_T$, independent of $y_0$ and $u$. 
\end{lemma}

\begin{remark} \rm Let $0 < T \le T_0$. By the arguments as in \Cref{rem-lem-transposition}, one can chose the constant $C_T$ in \eqref{lem-WP-Trans-cl} independent of $T$. 
\end{remark}

\begin{proof} By \Cref{pro-WP} and \Cref{lem-transposition}, it suffices to prove the uniqueness. Let $\mu > 0$ be large. We equip  $C^0([0, T]; \mH)$ with  the following norm 
$$
\vvvert  y \vvvert = \sup_{t \in [0, T]} e^{-\mu t} \| y(t) \|_{\mH}. 
$$
Recall that $y$ is a transposition solution if, for $t \ge 0$,  
\be\label{lem-WP-trans-F}
 \langle y (t), z_t \rangle - \langle y_0, z(0) \rangle =  \int_0^t \langle u (s), B^* z(s) \rangle \, ds 
+  \int_0^t \langle f(s),  z(s) \rangle \, ds + \int_0^t \langle y(s), M^* z (s) \rangle \, ds,
\ee
where $z_t \in \mH$ and $z$ is the weak solution of the backward system 
\be 
\left\{\begin{array}{c}
z' = - A^* z   \mbox{ in } (0, t), \\[6pt]
z(t) = z_t. 
\end{array} \right. 
\ee
Thus if $y$ and $\hy$ are two transposition solutions, then 
\be
 \langle y (t) - \hy (t), z_t \rangle  =  \int_0^t \langle y(s) - \hy(s), M^* z (s) \rangle \, ds.
\ee
This implies 
$$
e^{-\mu t}\|y(t) - \hy (t)\| \le C e^{-\mu t} \int_0^t \| y(s) - \hy(s)\| \, ds  \le \frac{C}{\mu} \vvvert y - \hy  \vvvert.  
$$
Here and in what follows in this proof, $C$ denotes a positive constant independent of $y$, $\hy$, and $\mu$. 
Thus
$$
\vvvert y - \hy \vvvert \le \frac{C}{\mu} \vvvert y -  \hy \vvvert. 
$$
The uniqueness follows and the proof is complete. 
\end{proof}

\begin{remark} \rm Similar results in the case $M=0$, $f = 0$, and $g = 0$ can be found in \cite[Section 2.3 of Chapter 2]{Coron07}. 

\end{remark}

\section{Dynamic feedback controls}\label{sect-dynamic-feedback}

This section consists of three subsections and is organized as follows. In the first subsection, we state and prove two useful lemmas, which will be used in the proofs of \Cref{thm-main1-D,thm-main1-D}. The proofs of \Cref{thm-main1-D} and \Cref{thm-main2-D} are given in the last two subsections, respectively.

\subsection{Two useful lemmas} Note that \eqref{identity-Op} can be written under an equivalent form as follows
\be \label{identity-Op-equiv}
A_\lambda Q + Q A_\lambda^* - B W B^* + QRQ = 0, 
\ee
where 
\be
A_\lambda = A + \lambda I.  
\ee
The meaning of \eqref{identity-Op-meaning} can be rewritten as follows 
\be \label{identity-Op-meaning-equiv} 
\langle Qx, A_\lambda^*y \rangle + \langle A_\lambda^*x, Q y \rangle  - \langle W B^*x,  B^*y \rangle  + \langle R Qx, Q y \rangle  = 0  \quad \forall \, x, y \in \cD(A^*). 
\ee

We have the following result concerning \eqref{identity-Op}. 

\begin{lemma}\label{lem-identity-Op} Assume \eqref{identity-Op}, i.e., \eqref{identity-Op-meaning}.  Given $x_0, y_0 \in \mH$ and $f, g \in L^1((0, T); \mH)$,  let $x, \, y \in C([0,T];\mH)$ be the unique weak solution of the systems 
$$
\left\{\begin{array}{c}
x' =  A_\lambda^*x + f \mbox{ in }  (0, T), \\[6pt]
x(0) = x_0,  
\end{array} \right. \quad \mbox{ and } \quad 
\left\{\begin{array}{c}
y' = A_\lambda^*y + g  \mbox{ in }  (0, T), \\[6pt]
y(0) = y_0.  
\end{array} \right.
$$
We have, for $t \in [0, T]$, 
\begin{multline} \label{lem-identity-Op-cl}
\langle Q x(t), y(t) \rangle -  \langle Q x_0, y_0 \rangle \\[6pt]= \int_0^t \Big( \langle W B^*x(s), B^*y (s) \rangle - \langle R Qx(s), Q y (s) \rangle \Big) \, d s + \int_0^t  \Big( \langle Q f(s), y(s) \rangle + \langle Q g(s), x(s) \rangle  \Big) \, ds. 
\end{multline}
\end{lemma}

\begin{proof} We first assume that $x_0, \,  y_0 \in \cD(A^*)$ and $f, \, g \in C([0, T]; \cD(A^*))$. 
Then $x, y \in C([0, T]; \cD(A^*))$ and $x', y' \in C([0, T]; \mH)$. We have
$$
\frac{d}{dt} \langle Q x, y \rangle = \langle  x', Q y \rangle + \langle Q x, y' \rangle = 
 \langle  A_\lambda^*x, Q y \rangle + \langle Q x, A_\lambda^*y  \rangle +   \langle  f, Q y \rangle + \langle Q x, g \rangle.    
$$
Using \eqref{identity-Op-meaning}, since $Q$ is symmetric, it follows that 
$$
\frac{d}{dt} \langle Q x, y \rangle =  \langle W B^*x, B^*y \rangle - \langle R Qx, Q y \rangle   + \langle Q f, y \rangle + \langle Q g, x \rangle.   
$$
We thus obtain \eqref{lem-identity-Op-cl}. 

\medskip
The proof in the general case is based on the previous case and a density argument using \Cref{lem-transposition}. 
\end{proof}

We next deal with the well-posedness of \eqref{thm-main1-sys} in \Cref{thm-main1-D}. It might be more convenient to consider a slightly more general system 
\be \label{CS-G-thm-main1-D}
\left\{\begin{array}{c}
y' = Ay + f  -  B W B^* \ty + M_1 y + M_2 \ty \mbox{ for } t \in (0, T), \\[6pt]
\ty' = - A^* \ty +\tf  + \tM_1 \ty + \tM_2 y  \mbox{ for } t \in (0, T), \\[6pt]
y(0) = y_0, \quad \ty(0) = \ty_0, 
\end{array} \right. 
\ee
with $y_0, \ty_0 \in \mH$,  $f, \tf \in L^1((0, T); \mH)$, $M_1, M_2, \tM_1, \tM_2 \in \cL(\mH)$, and $W \in \cL(\mU)$. As usual, a weak solution $(y, \ty)$ of \eqref{CS-G-thm-main1-D} is understood as an element $(y, \ty)\tr \in \big(C([0, T]; \mH) \big)^2$ such that 
\be\label{meaning-CS-G-thm-main1-D}
\left\{\begin{array}{c} \dsp \frac{d}{dt} \langle y, \varphi \rangle_{\mH}  = \langle A y + f  -  B W B^* \ty + M_1 y + M_2 \ty, \varphi \rangle_{\mH} \mbox{ in } [0, T] \\[6pt]
\dsp  \frac{d}{dt} \langle \ty, \tvarphi \rangle_{\mH}  = \langle -  A^* \ty + \tf + \tM_1 \ty + \tM_2 y, \tvarphi \rangle_{\mH} \mbox{ in } [0, T] \\[6pt]
\dsp y(0)  =  y_0,  \quad  \ty(0) = \ty_0,
\end{array}\right. \mbox{ for all } \varphi, \tvarphi \in \cD(A^*),  
\ee
for which 
\begin{itemize}
\item[$i)$] the differential equations are understood in the distributional sense, 

\item[$ii)$] the term $\langle A y + f  -  B W B^* \ty + M_1 y + M_2 \ty, \varphi \rangle_{\mH} $ is understood as $\langle y, A^*\varphi \rangle_{\mH}  +  \langle f + M_1 y + M_2 \ty, \varphi \rangle_{\mH} - \langle W B^* \ty, B^*\varphi \rangle_{\mU}$. 

\end{itemize}
Note that $B^* \ty \in L^2(0, T; \mU)$ since $B$ is an admissible control operator.   

\medskip 

We have the following result on the well-posedness of \eqref{CS-G-thm-main1-D}.

\begin{lemma}\label{lem-WP-thm-main1-D} Let $A$ be an infinitesimal generator of a {\it group}, and let $M_1, M_2, \tM_1, \tM_2 \in \cL(\mH)$ and $W \in \cL(\mU)$. Let $T>0$, $y_0, \ty_0 \in \mH$,  $f, \tf  \in L^1((0, T); \mH)$. There exists a unique weak solution $(y, \ty)\tr \in \big(C([0, T], \mH) \big)^2$ of \eqref{CS-G-thm-main1-D}. 
Moreover, with $g := f - B W B^* + M_1y + M_1 \ty$ and $\tg := \tf + \tM_1 \ty  + \tM_2 y$, we have \footnote{These identities below are understood in $\cD(A^*)'$ and $\cD(-A^*)'$, respectively.}
\be\label{pro-WP-identity-thm-main1-D-1}
y(t) = e^{t A} y_0 + \int_0^t e^{(t-s) A} g(s) \, ds \mbox{ for } t \in [0, T], 
\ee
and 
\be\label{pro-WP-identity-thm-main1-D-2}
\ty(t) = e^{- t A^*} y_0 + \int_0^t e^{-(t-s) A^*} \tg(s) \, ds \mbox{ for } t \in [0, T]. 
\ee
Moreover, we have 
$$
\| (y(t), \ty(t))\tr \|_{\mH} \le C \Big( \| (y_0, \ty_0)\tr \|_{\mH} + \| (f, \tf)\tr \|_{L^1((0, T); \mH)} \Big) \mbox{ in } [0, T], 
$$
for some positive constant $C$, independent of $y_0, \ty_0$, $f$, and $\tf$. 
\end{lemma}

\begin{remark} \rm In \Cref{lem-WP-thm-main1-D}, we does not require that $W$ is symmetric (or non-negative).  
\end{remark}

\begin{proof} We first note that $(y, \ty)\tr \in \big(C([0, T]; \mH) \big)^2$ is a weak solution of \eqref{CS-G-thm-main1-D} if and only if $(y, \ty)\tr \in \big(C([0, T]; \mH) \big)^2$, and \eqref{pro-WP-identity-thm-main1-D-1} and \eqref{pro-WP-identity-thm-main1-D-1} hold. This is a consequence of \Cref{pro-WP}.

We now establish the existence and uniqueness.  Let $\mu > 0$ be large. We equip  $\Big(C([0, T]; \mH) \Big)^2$ the following norm 
$$
\vvvert  y \vvvert  = \sup_{t \in [0, T]} e^{-\mu t} \| y(t) \|_{\mH}. 
$$

Define $\cF: \big(C([0, T]; \mH)\big)^2 \to \big(C([0, T]; \mH) \big)^2$ as follows 
\begin{equation*}
\cF 
\left(\begin{array}{c}
 y (t)  \\[6pt]
 \ty(t)
\end{array} \right) 
 = \left(\begin{array}{c}
\dsp e^{t A} y_0 + \int_0^t e^{(t-s) A} g(s) \, ds \\[6pt]
\dsp e^{-t A^*} \ty_0 + \int_0^t e^{- (t-s) A^*} \tg(s) \, ds 
\end{array}\right)
\mbox{ for } t \in [0, T]. 
\end{equation*}
Then, for $(y_1, \ty_1), (y_2, \ty_2) \in \big(C([0, T]; \mH)\big)^2$,  
\begin{multline*}
\cF \left(\begin{array}{c}
 y_2 (t)  \\[6pt]
 \ty_2(t)
\end{array} \right)  - \cF \left(\begin{array}{c}
 y_1 (t)  \\[6pt]
 \ty_1(t)
\end{array} \right)  \\[6pt]
=   \left(\begin{array}{c}
\dsp \int_0^t e^{(t-s) A} \Big(- B W B^* (\ty_2 - \ty_1) + M_1 (y_2 - y_1) + M_2 (\ty_2 - \ty_1)\Big) \, ds \\[6pt]
\dsp  \int_0^t e^{- (t-s) A^*} \Big( \tM_1(\ty_2 - \ty_1) + \tM_2 (y_2 - y_1) \Big) \, ds 
\end{array}\right). 
\end{multline*}
It follows that 
\begin{multline*}
 \left\|\cF \left(\begin{array}{c}
 y_2 (t)  \\[6pt]
 \ty_2(t)
\end{array} \right)  - \cF \left(\begin{array}{c}
 y_1 (t)  \\[6pt]
 \ty_1(t)
\end{array} \right) \right\|_{\mH} \\[6pt] 
\le C  \left( \int_0^t  \| (y_2, \ty_2)\tr(s) - (y_1, \ty_1)\tr(s)\|_{\mH}  \, ds + \|B^*(\ty_2 - \ty_1) \|_{L^2((0, t); \mU)} \right).
\end{multline*}
Here and in what follows in this proof, $C$ denotes a positive constant independent of solutions and $\mu$.  

This implies, by \eqref{lem-transposition-cl2} of \Cref{lem-transposition}, 
$$
e^{- \mu t} \left\|\cF \left(\begin{array}{c}
 y_2 (t)  \\[6pt]
 \ty_2(t)
\end{array} \right)  - \cF \left(\begin{array}{c}
 y_1 (t)  \\[6pt]
 \ty_1(t)
\end{array} \right) \right\|_{\mH} \\[6pt] 
\le C e^{-\mu t}  \int_0^t  \| (y_2, \ty_2)\tr(s) - (y_1, \ty_1)\tr (s)\|_{\mH}  \, ds. 
$$
We derive that 
$$
 \left\vvvert\cF \left(\begin{array}{c}
 y_2   \\[6pt]
 \ty_2
\end{array} \right)  - \cF \left(\begin{array}{c}
 y_1   \\[6pt]
 \ty_1
\end{array} \right) \right\vvvert \le \frac{C}{\mu}  \left\vvvert \left(\begin{array}{c}
 y_2   \\[6pt]
 \ty_2 
\end{array} \right)  - \left(\begin{array}{c}
 y_1   \\[6pt]
 \ty_1 
\end{array} \right) \right\vvvert.  
$$
By considering $\mu$ large enough, the existence and uniqueness of the weak solutions follow from a standard fixed point theorem. 
\end{proof}

\subsection{Dynamic feedback controls in the linear case - Proof of \Cref{thm-main1-D}} Set, for $t \ge 0$,  
\be
y_\lambda (t) = e^{\lambda t} y (t) \quad \mbox{ and } \quad \ty_\lambda (t) = e^{\lambda t} \ty (t), 
\ee
and denote
$$
A_\lambda = A + \lambda I. 
$$
We have
\be \label{thm-main1-sys-D}
\left\{\begin{array}{c}
y_\lambda' = A_\lambda y_\lambda - B W B^* \ty_\lambda \quad \mbox{ in } (0, + \infty), \\[6pt]
\ty_\lambda'  = - A_\lambda^* \ty_\lambda+ \lambda_1 Q^{-1} (y_\lambda - Q \ty_\lambda)   \quad \mbox{ in } (0, + \infty), \\[6pt]
y_\lambda(0) = y(0), \quad  \ty_\lambda(0) = \ty(0). 
\end{array} 
\right. 
\ee

Set, for $t  \ge 0$, 
$$
 Z_\lambda(t) = y_\lambda (t) - Q \ty_\lambda (t).
$$ 
We formally have, for $t \in (0, +\infty)$,   
\begin{multline*} 
\frac{d}{dt} Z_\lambda = A_\lambda y_\lambda - B W B^* \ty_\lambda + Q  A_\lambda^* \ty_\lambda  - \lambda_1 Z_{\lambda} \\[6pt]
= A_\lambda (y_\lambda  - Q \ty_\lambda) + A_\lambda Q  \ty_\lambda  - B W B^* \ty_\lambda+ Q  A_\lambda^* \ty_\lambda  - \lambda_1 Z_\lambda, 
\end{multline*}
which yields, since   \eqref{identity-Op} holds with  $R = 0$, that 
\be\label{thm-main1-p1-D}
\frac{d}{dt} Z_\lambda = A_\lambda Z_\lambda - \lambda_1 Z_\lambda. 
\ee

We now give the proof of \eqref{thm-main1-p1-D} (in the sense of weak solutions). 
Let  $\tau > 0$,  $\varphi_\tau \in \mH$ and let $\varphi \in C([0, \tau]; \mH)$ be the unique weak solution of the system 
\be \label{thm-main1-varphi-D} 
\left\{\begin{array}{c}
\varphi ' = - A_\lambda^* \varphi  \mbox{ in } (0, \tau), \\[6pt]
\varphi(\tau) = \varphi_\tau.  
\end{array} \right. 
\ee

Applying \Cref{lem-transposition} for $A_\lambda$ with $t = \tau$, we derive from \eqref{thm-main1-sys-D} and \eqref{thm-main1-varphi-D} that  
\be\label{thm-main1-p2-D}
 \langle y_\lambda(\tau), \varphi(\tau) \rangle - \langle y_\lambda (0), \varphi(0) \rangle =  - \int_0^\tau \langle W B^* \ty_\lambda (s), B^* \varphi(s) \rangle \, ds.
\ee
Applying \Cref{lem-identity-Op} for $A_\lambda$, $\ty_\lambda(\tau - \cdot)$, and $\varphi (\tau - \cdot)$ (with $R = 0$), we obtain   
\begin{multline}\label{thm-main1-p3-D}
\langle Q \ty_\lambda(0), \varphi(0) \rangle - \langle Q \ty_\lambda(\tau), \varphi (\tau) \rangle  \\[6pt]
 =  \int_0^\tau  \langle W B^*\ty_\lambda(\tau - s), B^*\varphi (\tau - s) \rangle  \, d s - \lambda_1 \int_0^\tau \langle Z_\lambda(\tau - s),  \varphi(\tau - s) \rangle  \, ds.   
\end{multline}
Summing \eqref{thm-main1-p2-D} and \eqref{thm-main1-p3-D},   we deduce from \eqref{thm-main1-sys-D} and \eqref{thm-main1-varphi-D} that  
$$
\langle Z_\lambda(\tau), \varphi(\tau) \rangle - \langle Z_\lambda(0), \varphi(0) \rangle = - \lambda_1 \int_0^\tau \langle Z_\lambda(\tau - s), \varphi(\tau - s) \rangle  \, ds. 
$$
This yields 
$$
\langle Z_\lambda(\tau), \varphi(\tau) \rangle - \langle Z_\lambda(0), e^{\tau A^*}\varphi(\tau) \rangle = - \lambda_1 \int_0^\tau \langle Z_\lambda(\tau - s), e^{s A^*}\varphi(\tau) \rangle  \, ds. 
$$
Since $\varphi(\tau) \in \mH$ is arbitrary, we obtain
$$
Z_\lambda(\tau) = e^{\tau A} Z_\lambda(0) - \lambda_1  \int_0^\tau e^{(\tau - s) A} Z_\lambda(s) \, ds, 
$$
which implies \eqref{thm-main1-p1-D}.

We derive from \eqref{thm-main1-p1-D} that 
\be \label{thm-main1-p00-D}
\| Z_\lambda (t) \|_{\mH} \le C e^{(-\lambda_1 + \lambda + \homega_0(A))t} \|Z_\lambda(0) \|_{\mH},    
\ee
which yields 
\be \label{thm-main1-p11-D}
\| y(t) - Q \ty(t) \|_{\mH} \le C e^{(-\lambda_1 + \homega_0(A))t} \| y(0) - Q \ty(0) \|_{\mH}.  
\ee
Here and in what follows in this proof, $C$ is a positive constant independent of $t$ and $(y_0, \ty_0)$. 

Since 
$$
\ty'  = - A^* \ty - 2 \lambda \ty + \lambda_1 Q^{-1} (y - Q\ty) \quad \mbox{ in } (0, + \infty), 
$$
it follows that 
$$
\ty_{2 \lambda}' = - A^* \ty_{2 \lambda} + f (t) \mbox{ in } (0, + \infty), 
$$
where
$$
\ty_{2 \lambda} = e^{2 \lambda} \ty (t) \quad \mbox{ and } \quad f(t) = \lambda_1 e^{2 \lambda t} Q^{-1} (y(t) - Q \ty (t)) \mbox{ in } (0, + \infty).  
$$
We obtain 
\be\label{thm-main1-p4-D}
\ty_{2 \lambda}(t) = e^{- t A^*} \ty_{2 \lambda}(0) + \int_0^t e^{-(t-s) A^* }f(s) \, ds. 
\ee

From the definition of $f$ and \eqref{thm-main1-p11-D}, we have 
$$
\left\| \int_0^t e^{-(t-s) A^* }f(s) \, ds \right \|_{\mH} \le C \int_0^t e^{\homega_0 (-A^*) (t-s) } e^{\big(-\lambda_1 + \homega_0(A) + 2 \lambda \big) s} \|y(0) - Q \ty(0) \|_{\mH} \, ds.  
$$
Since
$$
-\homega_0(-A^*)  + \homega_0(A) + 2 \lambda - \lambda_1 \mathop{<}^{\eqref{thm-main1-lambda-D}}  0,  
$$
it follows that 
\be\label{thm-main1-p5-D}
\left \| \int_0^t e^{-(t-s) A^* }f(s) \, ds \right \|_{\mH} 
\le C e^{ \homega_0 (-A^*) t} \|y(0) - Q \ty(0) \|_{\mH}. 
\ee

Combining \eqref{thm-main1-p4-D} and \eqref{thm-main1-p5-D} yields 
$$
\| \ty_{2 \lambda} (t) \|_{\mH} \le C e^{\homega_0(-A^*) t} \big( \| \ty(0) \|_{\mH} + \|y(0) - Q \ty_0 \|_{\mH} \big), 
$$
which implies 
\be\label{thm-main1-p22-D}
\| \ty(t)\|_{\mH} \le C e^{(\homega_0(-A^*) - 2 \lambda) t} \big( \| y(0) \|_{\mH} + \|\ty(0) \|_{\mH} \big).  
\ee

Combining \eqref{thm-main1-p11-D} and \eqref{thm-main1-p22-D}, we obtain 
\be\label{thm-main1-p33-D}
\| y(t)\|_{\mH} + \| \ty(t)\|_{\mH} \le C\Big(  e^{(-\lambda_1 + \homega_0(A))t}  + e^{(\homega_0(-A^*) - 2 \lambda) t}  \Big) \big( \| y(0) \|_{\mH} + \|\ty(0)\|_{\mH} \big).  
\ee
Since 
$$
\lambda_1 - \homega_0(A) \mathop{>}^{\eqref{thm-main1-lambda-D}} 2 \lambda - \homega_0(-A^*),  
$$
it follows from \eqref{thm-main1-p33-D} that 
\be
\| y(t)\|_{\mH} + \| \ty(t)\|_{\mH} \le C  e^{(\homega_0(-A^*) - 2 \lambda) t}  \big( \| y(0) \|_{\mH} + \|\ty(0)\|_{\mH} \big), 
\ee
which is \eqref{thm-main1-cl1-D}. 

\medskip
It is clear that \eqref{thm-main1-cl2-D} is a direct consequence of \eqref{thm-main1-cl1-D}. 

\medskip
The proof is complete. 
 \qed

\subsection{Dynamic feedback controls in the nonlinear case - Proof of \Cref{thm-main2-D}} 

For each $T > 0$, there exists $\eps = \eps_T > 0$ such that  \eqref{thm-main2-sys-D} is well-posed in the time interval $[0, T]$. The global existence and uniqueness follow for small $\eps$  provided that 
\eqref{thm-main2-cl1-D} is established for each fixed time interval $[0, T]$ with $\eps_T$ sufficiently small. The proof is in the same spirit of the one of \Cref{thm-main1-D} but more involved due to the nonlinearity. 

Set, for $t \ge 0$,  
\be
y_\lambda (t) = e^{\lambda t} y (t) \quad \mbox{ and } \quad \ty_\lambda (t) = e^{\lambda t} \ty (t), 
\ee
and denote
$$
A_\lambda = A + \lambda I. 
$$
We have
\be 
\left\{\begin{array}{c}
y_\lambda' = A_\lambda y_\lambda + e^{\lambda \cdot} f(e^{-\lambda \cdot} y_\lambda )- B W B^* \ty_\lambda \quad \mbox{ in } (0, + \infty), \\[6pt]
\ty_\lambda'  = - A_\lambda^* \ty_\lambda + Q^{-1} e^{\lambda \cdot} f (e^{-\lambda \cdot} Q \ty_{\lambda} )+ \lambda_1 Q^{-1} (y_\lambda - Q \ty_\lambda)   \quad \mbox{ in } (0, + \infty), \\[6pt]
y_\lambda(0) = y(0), \quad  \ty_\lambda(0) = \ty(0). 
\end{array} 
\right. 
\ee

Set, for $t  \ge 0$, 
$$
 Z_\lambda(t) = y_\lambda (t) - Q \ty_\lambda (t).
$$ 
As in the proof of \eqref{thm-main1-p1-D} in the proof of \Cref{thm-main1-D}, we derive that $Z_\lambda$ is a weak solution of the equation 
\be \label{thm-main2-p1-D}
\frac{d}{dt} Z_\lambda =   A_\lambda Z_\lambda - \lambda_1 Z_\lambda + g_1, 
\ee
where 
$$
g_1(t) = e^{\lambda t} \Big( f(e^{-\lambda t} y_\lambda (t) ) - f (e^{-\lambda t} Q \ty_{\lambda}(t) ) \Big) \mbox{ for } t \in (0, +\infty). 
$$

It follows from \eqref{thm-main2-p1-D} that 
\be \label{thm-main2-p00-D}
\| Z_\lambda (t) \|_{\mH} \le C e^{(-\lambda_1 + \lambda + \homega_0(A))t} \|Z_\lambda(0) \|_{\mH} + C \int_0^t e^{(-\lambda_1 + \lambda + \homega_0(A))(t-s)} \| g_1(s)\|_{\mH} \, ds.   
\ee
Here and in what follows in this proof, $C$ is a positive constant independent of $t$ and $(y_0, \ty_0) \tr$.

From \eqref{thm-main2-p00-D}, we obtain 
\begin{multline} \label{thm-main2-p11-D}
\| y(t) - Q \ty(t) \|_{\mH} \\[6pt]
\le C e^{(-\lambda_1 + \homega_0(A))t} \| y(0) - Q \ty(0) \|_{\mH} + C e^{(-\lambda_1 + \homega_0(A))t}  \int_0^t e^{-(-\lambda_1 + \lambda + \homega_0(A))s} \| g_1(s)\|_{\mH} \, ds.  
\end{multline}

Since 
$$
\ty'  = - A^* \ty - 2 \lambda \ty + e^{\lambda \cdot} Q^{-1}  f (e^{-\lambda \cdot} Q \ty_{\lambda} ) + \lambda_1 Q^{-1} (y - Q\ty) \quad \mbox{ in } (0, + \infty), 
$$
it follows that 
$$
\ty_{2 \lambda}' = - A^* \ty_{2 \lambda} + f_1(t) + f (t) \mbox{ in } (0, + \infty), 
$$
where, in $(0, + \infty)$, 
$$
\ty_{2 \lambda} = e^{2 \lambda} \ty (t),  \quad h(t) = \lambda_1 e^{2 \lambda t} Q^{-1} (y(t) - Q \ty (t)), \quad \mbox{ and } \quad h_1(t) =  e^{3\lambda t}  Q^{-1} f (e^{-\lambda t} Q \ty_{\lambda} (t)). 
$$
We derive that 
\be\label{thm-main2-p4-D}
\ty_{2 \lambda}(t) = e^{- t A^*} \ty_{2 \lambda}(0) + \int_0^t e^{-(t-s) A^* }\big(h(s) + h_1(s) \big) \, ds. 
\ee

Using the first inequality in \eqref{thm-main2-lambda-D},  we derive from \eqref{thm-main2-p11-D} that 
\begin{multline}\label{thm-main2-p5-D}
\left \| \int_0^t e^{-(t-s) A^* } h(s) \, ds \right \|_{\mH} \\[6pt]
\le C e^{t \homega_0 (-A^*)} \Big( \|y(0) - Q \ty(0) \|_{\mH} + \int_0^t e^{-(-\lambda_1 + \lambda + \homega_0(A))s} \| g_1(s)\|_{\mH} \, ds \Big) \mbox{ for } t \ge 0.  
\end{multline}
Using \eqref{cond-NL1} and the first inequality in \eqref{thm-main2-lambda-D}, we derive from \eqref{thm-main2-p11-D}, \eqref{thm-main2-p4-D}, and \eqref{thm-main2-p5-D} that for every $\eps > 0$, there exists $\delta > 0$ such that if $\| (y(t), \ty(t)) \|_{\mH} \le \delta$ in $[0, T]$ for some $T>0$, then 
\begin{multline*}
\|(y(t), \ty (t)) \|_{\mH} \le C e^{ (\homega_0(-A^*) - 2 \lambda) t} \|(y_0, \ty_0) \|_{\mH} 
\\[6pt] 
+ C \eps e^{ (\homega_0(-A^*) - 2 \lambda) t} \int_0^t e^{-(-\lambda_1 + \lambda + \homega_0(A))s} \|(y(s), \ty (s)) \|_{\mH} \, ds \\[6pt]
+ C \eps e^{ (\homega_0(-A^*) - 2 \lambda) t}  \int_0^t e^{\big(3 \lambda - \homega_0(-A^*)\big) s} \|(y(s), \ty (s)) \|_{\mH} \, ds \mbox{ for } t \in [0, T]. 
\end{multline*}
Here $C$ is a positive constant independent of $T$, $\eps$ and $\delta$. Thus, 
for all $T > 0$, there exists $\delta > 0$ such that if $\| (y_0, \ty_0) \|_{\mH} \le \delta$ then
\be\label{thm-main2-D-cc1}
\|(y(t), \ty (t)) \|_{\mH} \le C e^{ (\homega_0(-A^*) - 2 \lambda) t} \|(y_0, \ty_0) \|_{\mH}  \mbox{ in } [0, T]. 
\ee
In particular, we derive that if $T$ is chosen sufficiently large,  
\be\label{thm-main2-D-cc2}
\|(y(T), \ty (T)) \|_{\mH} \le e^{ (\omega_0(-A^*) - 2 \gamma) T} \|(y_0, \ty_0) \|_{\mH}. 
\ee
The conclusion follows from \eqref{thm-main2-D-cc1} and \eqref{thm-main2-D-cc2} by considering the time $n T \le t \le n (T+1) $ for $n \in \N$. \qed

\section{Static feedback controls in the trajectory sense} \label{sect-static-feedback}
 
This section consisting of three subsections is organized as follows. The proofs of \Cref{thm-main1} and \Cref{thm-main2} are given in the first two subsections, respectively. In the last subsection, we study of the infinitesimal generator of the semigroup associated with the static feedback controls given in \Cref{thm-main1}.

\subsection{Static feedback controls in the linear case - Proof of \Cref{thm-main1}} Set, for $t \ge 0$,  
\be
y_\lambda (t) = e^{\lambda t} y (t) \quad \mbox{ and } \quad \ty_\lambda (t) = e^{\lambda t} \ty (t).
\ee
We then have, with $A_\lambda = A + \lambda I$,  
\be 
\left\{\begin{array}{c}
y_\lambda' = A_\lambda y_\lambda - B W B^* \ty_\lambda \quad \mbox{ in } (0, + \infty), \\[6pt]
\ty_\lambda'  = - A_\lambda^* \ty_\lambda - R Q \ty_\lambda \quad \mbox{ in } (0, + \infty), \\[6pt]
y_\lambda(0) = y(0), \quad  \ty_\lambda(0) = \ty(0) (= Q^{-1} y(0)). 
\end{array} 
\right. 
\ee

Set, for $t  \ge 0$, 
$$
 Z_\lambda(t) = y_\lambda (t) - \tz_\lambda (t).
$$ 
As in the proof of \eqref{thm-main1-p1-D} in the proof of \Cref{thm-main1-D}, we derive that $Z_
\lambda$ is a weak solution of the equation  
\be\label{thm-main1-Z}
\frac{d}{dt} Z_\lambda   = A_\lambda Z_\lambda. 
\ee
Since $Z_\lambda(0) = 0$, it follows that 
\be \label{thm-main1-p1}
Z_\lambda (t) = 0 \mbox{ for } t \ge 0.  
\ee
In other words, \eqref{thm-main1-cl1} holds.

We next deal with \eqref{thm-main1-cl2}. Formally, we have
\begin{multline}\label{thm-main1-cccc}
\frac{d}{dt} \langle y, \ty \rangle =  \langle A y - B W B^* \ty, \ty \rangle_{\mH} + \langle y, - A^*\ty - R Q \ty - 2 \lambda \ty \rangle_{\mH}  \\[6pt]  \mathop{=}^{\eqref{thm-main1-cl1}} 
  \langle A y - B W B^* \ty, \ty \rangle_{\mH} + \langle y, - A^*\ty - R y - 2 \lambda Q^{-1} y \rangle_{\mH} \\[6pt]
  = -  \| W^{1/2} B^*\ty \|_{\mU}^2 - \|R^{1/2} y\|_{\mH}^2 - 2 \lambda \langle Q^{-1} y, y \rangle_{\mH},  
\end{multline}
which yields \eqref{thm-main1-cl2}. The rigor proof of \eqref{thm-main1-cl2} can be done by applying \Cref{lem-transposition} for $y$ and $\ty$. 

To derive \eqref{thm-main1-cl3} from \eqref{thm-main1-cl2}, one just needs to set 
$$
\rho(t) =   \langle Q^{-1} y(t), y(t) \rangle \mbox{ for } t \ge 0, 
$$
and note that, by \eqref{thm-main1-cl2}, 
$$
\rho \in W^{1,1}(0, T) \mbox{ for all } T > 0 \quad \mbox{ and } \quad  \rho'(t) \le - 2 \lambda \rho (t) \mbox{ for } t \ge 0.  
$$

The proof is complete. \qed

\subsection{Static feedback controls in the nonlinear case - Proof of \Cref{thm-main2}} For each $T > 0$ there exists $\eps = \eps_T > 0$ such that  \eqref{thm-main2-sys} is well-posed in the time interval $[0, T]$. The global existence and uniqueness follow for small $\eps$  provided that 
\eqref{thm-main2-cl1}, \eqref{thm-main2-cl2}, and \eqref{thm-main2-cl3} are established for each fixed time interval $[0, T]$ with $\eps_T$ sufficiently small.  

We now establish \eqref{thm-main2-cl1}, \eqref{thm-main2-cl2}, \eqref{thm-main2-cl3} in $[0, T]$ for $\eps < \eps_T$ (small). Set, in $[0, T]$,  
\be
y_\lambda (t) = e^{\lambda t} y (t) \quad \mbox{ and } \quad \ty_\lambda (t) = e^{\lambda t} \ty (t).
\ee
We then have, with $A_\lambda = A + \lambda I$,  
\be \label{thm-main2-sys}
\left\{\begin{array}{c}
y_\lambda' = A_\lambda y_\lambda + e^{\lambda t} f(e^{-\lambda t} y_\lambda (t)) - B W B^* \ty_\lambda \quad \mbox{ in } (0, T), \\[6pt]
\ty_\lambda'  = - A_\lambda^* \ty_\lambda - R Q \ty_\lambda +  Q^{-1} e^{\lambda t} f (e^{-\lambda t} Q \ty_\lambda) \quad \mbox{ in } (0, T), \\[6pt]
y_\lambda(0) = y(0), \quad  \ty_\lambda(0) = \ty(0) (= Q^{-1} y(0)). 
\end{array} 
\right. 
\ee

Set, for $t \in [0, T]$, 
$$
Z_\lambda(t) = y_\lambda (t) - Q \ty_\lambda (t).
$$ 
As in the proof of \eqref{thm-main1-p1-D} in the proof of \Cref{thm-main1-D}, we derive that  
\be\label{thm-main2-Z}
\frac{d}{dt} Z_\lambda =  A_\lambda Z_\lambda + e^{\lambda \cdot} \Big(f(e^{-\lambda \cdot} y_\lambda ) - f(e^{-\lambda \cdot} Q \ty_\lambda (\cdot)) \Big). 
\ee
Since $Z_\lambda (0) = 0$, we obtain 
$$
y_\lambda(t) - Q \ty_\lambda(t) = \int_0^t  e^{(t-s) A_\lambda }e^{\lambda s} \Big(f(e^{-\lambda s} y_\lambda (s)) - f(e^{-\lambda s} Q \ty_\lambda (s)) \Big) \, ds. 
$$
Using \eqref{cond-NL2}, we deduce that 
$$
y_\lambda(t) = Q \ty_\lambda(t) \mbox{ for } t \ge 0,
$$
which implies \eqref{thm-main2-cl1}. 

\medskip
We next deal with \eqref{thm-main2-cl2}. The proof of \eqref{thm-main2-cl2} is similar to the one of \eqref{thm-main1-cl2} by applying \Cref{lem-transposition} for $y$ and $\ty$. 

\medskip
What have been done so far does not require $\lambda > 0$. The fact $\lambda > 0$ is used to derive \eqref{thm-main2-cl3} from \eqref{thm-main2-cl2}. Set 
$$
\rho(t) =   \langle Q^{-1} y(t), y(t) \rangle \mbox{ for } t \ge 0, 
$$
Note that, by \eqref{thm-main2-cl2}, as in the proof of \eqref{thm-main2-D-cc1} for all $T> 0$, there exists $\delta > 0$ such that if $\| y_0 \|_{\mH} \le \delta $ in $[0, T]$, then 
\be \label{thm-main2-cc1}
\rho(t) \le C e^{- 2 \lambda t} \rho(0) \mbox{ in } [0, T].   
\ee
In particular, we have, if $T$ is chosen sufficiently large,  
\be\label{thm-main2-cc2}
\|(y(T), \ty (T)) \|_{\mH} \le e^{- 2 \gamma T} \|(y_0, \ty_0) \|_{\mH}. 
\ee
The conclusion follows from \eqref{thm-main2-cc1} and \eqref{thm-main2-cc2} by considering the time $n T \le t \le n (T+1) $ for $n \in \N$. 

\medskip
The proof is complete. 
\qed

\subsection{The infinitesimal generator of the semigroup associated with the static feedback controls} \label{sect-generator}

Here is the main result of this section on the infinitesimal generator of the semigroup associated with the static feedback controls from \Cref{thm-main1}

\begin{proposition} \label{pro-main1} Let $\lambda \in \mR$ and assume \eqref{identity-Op}. Let $y_0 \in \mH$, set 
\be \label{thm-main1-semigroup}
S^Q(t)(y_0) = y(t), 
\ee
where $(y, \ty)$ is the solution of \eqref{thm-main1-sys}. Then 
\be \label{pro-main1-SG1}
\big(S^Q(t) \big)_{t \ge 0} \mbox{ is a strongly continuous semigroup on $\mH$}. 
\ee
Moreover, the semigroup $\big(S^Q(t) \big)_{t \ge 0}$ decays exponentially with the rate $\lambda$, i.e., there exists $C>0$ such that 
\be \label{pro-main1-SG2}
\|S^Q(t)\| \le C e^{-\lambda t} \mbox{ for } t \ge 0. 
\ee
Let $(A^Q, \cD(A^Q))$ be its infinitesimal generator.  We have 
\be\label{pro-main1-cl1}
\cD(A^Q) = Q \cD(A^*) : =  \Big\{Q x ; x \in \cD(A^*) \Big\}
\ee
and 
\be \label{pro-main1-cl2}
A^Q z = -Q A^* Q^{-1} z - 2 \lambda z - Q R z \mbox{ for } z \in \cD(A^Q).  
\ee
We also have
\begin{itemize}
\item[$i)$] if $BWB^*$ is bounded, i.e.,  $BWB^* \in \cL(\mH)$, then 
\be \label{pro-main1-cl3}
\cD(A^Q) = \cD(A) \quad \mbox{ and } \quad  A^Q x = A x - B W B^* Q^{-1} x \mbox{ for } x \in \cD(A) = \cD(A^Q). 
\ee

\item[$ii)$] if $\cD(A^Q) = \cD(A)$, then $BWB^*x \in \mH$ for $x \in \cD(A^*)$,  and 
\be \label{pro-main1-cl4}
\| B W B^*x \|_{\mH} \le \| A Q x\|_{\mH}  + C ( \| A^*x \|_{\mH} + \| x \|_{\mH}) \mbox{ for } x \in \cD(A^*)
\ee
for some positive constant $C$ independent of $x$. 
\end{itemize}
\end{proposition}

\begin{proof}[Proof of \Cref{pro-main1}] 

It is clear that \eqref{pro-main1-SG1} and \eqref{pro-main1-SG2} are the consequences of \Cref{thm-main1}. 

\medskip 
We now prove \eqref{pro-main1-cl1} and \eqref{pro-main1-cl2}. 
Fix $y_0 \in Q \cD(A^*)$ (arbitrary). Let $(y, \ty)$ be the unique weak solution of \eqref{thm-main1-sys}. Since $\ty(0) = Q^{-1} y_0 \in \cD(A^*)$, it follows that $\ty \in C^1([0, + \infty); \mH) \cap C^0([0, + \infty); \cD(A^*))$ and  
\be \label{pro-main1-p1}
\ty'(0) =  A^* \ty(0) - 2 \lambda \ty(0) - R Q \ty(0). 
\ee
Since $y(t) = Q \ty (t)$ for $t \ge 0$ by \Cref{thm-main1}, we derive that $y'(0)$ is well-defined and 
$$
y'(0) = Q \ty'(0) \mathop{=}^{\eqref{pro-main1-p1}} -Q A^* Q^{-1} y_0 - 2 \lambda y_0 - Q R Q y_0. 
$$ 
Hence $y_0 \in \cD(A^Q)$ and 
$$
A^Q y_0 = - Q A^* Q^{-1} y_0 - 2 \lambda y_0 - Q R Q y_0. 
$$
To complete the proof of  \eqref{pro-main1-cl1} and \eqref{pro-main1-cl2}, we now show that if $y_0 \in \cD(A^Q)$ then $y_0 \in Q \cD(A^*)$. Fix $y_0 \in \cD(A^Q)$ (arbitrary) and let  $(y, \ty)$ be the unique solution of \eqref{thm-main1-sys}. Since $y_0 \in \cD(A^Q)$ and $S^Q(t)(y_0) = y(t)$, it follows that $y \in  C^1([0, + \infty); \mH) \cap C^0([0, + \infty); \cD(A^Q))$. In particular $y'(0)$ is well-defined. Since $y(t)  = Q \ty(t)$ for $t \ge 0$ by \Cref{thm-main1}, it follows from the equation of $\ty$ in \eqref{thm-main1-sys} that $\ty'(0)$ is well-defined and thus $\ty(0) \in \cD(A^*)$. Since $\ty(0) = Q^{-1} y_0$, we derive that 
$$
Q^{-1} y_0  \in \cD(A^*). 
$$
In other words, $y_0 \in Q \cD(A^*)$. 

\medskip

We next establish \eqref{pro-main1-cl3}. We first assume that $B W B^* \in \cL(\mH)$.  It follows that the generator of the semigroup $\big(S^Q(t) \big)_{t \ge 0}$ is $A - B W B^* Q^{-1}$ with the domain $\cD(A)$. 

\medskip
We finally derive \eqref{pro-main1-cl4}. Assume that $\cD(A^Q) = \cD(A)$.  From \eqref{identity-Op-meaning}, we have, for $x, y \in \cD(A^*)$, 
\begin{align*}
|\langle W B^*x,  B^*y \rangle_{\mU}|  \le &  |\langle Qx, A^*y \rangle_{\mH}| +  | \langle Q y, A^*x \rangle_{\mH}|  + |\langle R Qx, Q y \rangle_{\mH}| + 2 |\lambda| | \langle Q x, y \rangle_{\mH}| \\[6pt]
\le &  (\| A Q x\|_{\mH}  + C \| A^*x \|_{\mH} + C \| x \|_{\mH})   \| y \|_{\mH}. 
\end{align*}
It follows that 
$$
\| B W B^*x \|_{\mH} \le \| A Q x\|_{\mH}  + C ( \| A^*x \|_{\mH} + \| x \|_{\mH}) \mbox{ for } x \in \cD(A^*), 
$$
which is \eqref{pro-main1-cl4}. 

\medskip 
The proof is complete. 
\end{proof}

\begin{remark} \rm Related results to \Cref{pro-main1} from the linear quadratic optimal control theory can be found in  \cite{Flandoli87-LQR, CZ95, BLT00, WR00,Triggiani05}.  Known results established in the case  $\lambda =0$ and $W$ being identity are connections between $\cD({A^Q}^*)$ and $\cD(A)$, see \cite[Theorem 2.1]{Triggiani05}. This is different from \eqref{pro-main1-cl1} where a connection between $\cD(A^Q)$ and $\cD(A^*)$ is established. Assertion $i)$ is equivalent to the fact that $B$ is bounded, i.e., $B \in \cL(\mU, \mH)$ when $W$ is positive; this case is well-known. 
\end{remark}

\section{Choices of $Q$ for exactly controllable systems} \label{sect-Q-E}

In this section, we discuss how to choose $Q$ for exactly controllable systems. Assume that the system is exactly controllable at time $T$. This is equivalent to the fact that \eqref{observability-inequality-CS} holds.  Fix $\lambda \in \mR$ and  $T_* > T$ and let $\rho: [0, T_*] \to \mR$ be such that 
\be \label{rho-E}
\mbox{$\rho$ is Lipschitz, $\rho$ is decreasing, $\rho (0) = 1$, $\rho(T) > 0$,  and $\rho (T_*) = 0$.} 
\ee
Let $W  \in \cL(\mU)$ be symmetric and positive.  Define $Q: \mH \to \mH$ as follows   
\be \label{def-Q-E}
\langle Q z_1, z_2 \rangle= \int_{0}^{T_*} \rho(s) e^{-2 \lambda s} \langle W B^* e^{-s A^*} z_1, B^* e^{-s A^*} z_2 \rangle \, ds \mbox{ for } z_1, z_2 \in \mH.   
\ee
Then $Q$ is linear, continuous, and symmetric. Moreover, since $\rho$ is decreasing and $\rho(T) > 0$, $A$ is an infinitesimal of a group,  it follows from \eqref{observability-inequality-CS} that 
\be
Q \mbox{ is invertible}. 
\ee 
Let $R: \mH \to \mH$ be defined by 
\be\label{def-N-E} 
\langle  R Qz_1, Q z_2 \rangle= - \int_0^{T_*} \rho'(s) \langle W B^* e^{-s (A + \lambda I)^*} z_1, B^* e^{-s (A + \lambda I)^*} z_2 \rangle \, ds. 
\ee

For $z_1, z_2 \in \cD({A^*}^2)$, we have, from \eqref{def-Q-E},  
\begin{multline}\label{E-p1}
\langle Q z_1, (A + \lambda I)^* z_2 \rangle + \langle (A + \lambda I)^* z_1, Q z_2 \rangle \\[6pt]
= \int_{0}^{T_*} \rho(s)  \langle W B^* e^{-s (A + \lambda I)^*} z_1, B^* e^{-s (A + \lambda I)^*} (A + \lambda I)^* z_2 \rangle \, ds \\[6pt]
+  \int_{0}^{T_*} \rho(s)  \langle W B^* e^{-s (A+ \lambda I)^*} (A + \lambda I)^* z_1, B^* e^{-s (A + \lambda I)^*}  z_2 \rangle \, ds. 
\end{multline}
Using the fact that, for $z \in \cD({A^*}^2)$, 
$$
e^{-s (A+ \lambda I)^*} (A + \lambda I)^* z = - \frac{d}{ds} \left( e^{-s (A+ \lambda I)^*} z \right), 
$$ 
we derive from \eqref{E-p1} that 
\begin{multline}
\langle Q z_1, (A + \lambda I)^* z_2 \rangle + \langle (A + \lambda I)^*Q z_1, z_2 \rangle
= -  \int_{0}^{T_*} \rho(s)  \frac{d}{ds} \Big( e^{-s (A + \lambda I)} B W B^* e^{-s (A + \lambda I)^*} \Big) \, ds, 
\end{multline}
which yields, by an integration by parts,  
\begin{multline}
\langle Q z_1, (A + \lambda I)^* z_2 \rangle + \langle (A + \lambda I)^*Q z_1, z_2 \rangle
 \\[6pt]
=   \langle W B^* z_1, B^* z_2 \rangle +  \int_0^{T_*} \rho'(s) \langle W B^* e^{-s (A + \lambda I)^*} z_1, B^* e^{-s (A + \lambda I)^*} z_2 \rangle \, ds, 
\end{multline}
This implies \eqref{identity-Op-meaning} for $z_1, z_2 \in \cD({A^*}^2)$. The general case follows by density. 

\medskip

We have just proven the following result. 

\begin{proposition} \label{pro-Q-E} Assume that $(S(t))_{t \in \mR} \subset \cL(\mH)$  is a strongly continuous group in $\mH$, $B$ is an admissible control operator, and system \eqref{CS} is exactly controllable in time $T$ for some $T>0$. Let $\lambda \in \mR$, $T_* > T$, and $\rho: [0, T_*] \to \mR$ be a function satisfying \eqref{rho-E}, and let $W \in \cL(\mU)$ be symmetric and positive. Define $Q: \mH \to \mH$ by 
\be 
\langle Q z_1, z_2 \rangle= \int_{0}^{T_*} \rho(s) e^{-2 \lambda s} \langle W B^* e^{-s A^*} z_1, B^* e^{-s A^*} z_2 \rangle \, ds \mbox{ for } z_1, z_2 \in \mH.   
\ee
Then $Q$ is linear, continuous, symmetric, and invertible and  \eqref{identity-Op} holds with $R$ being defined by \eqref{def-N-E}, i.e., \eqref{identity-Op-meaning} is valid. 
\end{proposition}

\begin{remark} \label{rem-Q-E} \rm \Cref{pro-Q-E} covers the setting considered by Komornik. Indeed, set, with $T_* = T + \frac{1}{2 \lambda}$ 
\be \label{def-rho-K}
\rho(t) = \left\{ \begin{array}{cl}
1 & \mbox{ for } 0 \le t \le T, \\[6pt]
2 \lambda e^{-2 \lambda (T - t)} (T_* - t) &  \mbox{ for } T <  t \le T_*. 
\end{array}\right.
\ee
Then 
$$
e_\lambda(t) = e^{\lambda t} \rho(t) \mbox{ in } [0, T_*]. 
$$
Since, for $T \le t \le T_* = T + \frac{1}{2 \lambda}$,  
$$
\rho (t) = e \tau e^{-\tau} \mbox{ with } \tau = 2 \lambda (T_* - t), 
$$
and the function $\tau e^{-\tau}$ is increasing in $[0, 1]$, it follows that $\rho$ defined in \eqref{def-rho-K} verifies \eqref{rho-E}.
\end{remark}

When $A$ is skew-adjoint and $R = 0$, one has the following result.   

\begin{proposition} \label{pro-Q-E-U} Assume that $\big(S(t) \big)_{t \in \mR} \subset \cL(\mH)$ is  a strongly continuous group, $B$ is an admissible control operator, and system \eqref{CS} is exactly controllable in time $T$ for some $T>0$. Let $\lambda \in \mR$ and let $W \in \cL(\mU)$ be symmetric and non-negative, and assume that $\lambda > \omega_0(-A^*)$. 
Define $Q: \mH \to \mH$ by 
\be \label{pro-Q-E-U-defQ}
\langle Q z_1, z_2 \rangle= \int_{0}^{\infty} e^{-2 \lambda s} \langle W B^* e^{-s A^*} z_1, B^* e^{-s A^*} z_2 \rangle \, ds \mbox{ for } z_1, z_2 \in \mH.   
\ee
Then $Q$ is linear, continuous, symmetric, and invertible, and \eqref{identity-Op} holds with $R = 0 $, i.e., \eqref{identity-Op-meaning} is valid with $R = 0$. 
\end{proposition}

\begin{proof} The proof of \eqref{pro-Q-E} is almost the same as the one of \Cref{pro-Q-E}. One just needs to note that the RHS of \eqref{pro-Q-E-U-defQ} is well-defined for $\lambda > \omega_0(-A^*)$. The details are omitted. 
\end{proof}

\begin{remark} \rm
\Cref{pro-Q-E-U} was previously obtained by Urquiza \cite{Urquiza05} by a different approach using results of Grabowski in \cite{Grabowski90} (see also \cite{HW97}).
\end{remark}

\noindent {\bf Acknowledgement.} The author would like to thank Jean-Michel Coron, Emmanuel Tr\'elat, and Marius Tusnack for many interesting discussions, useful comments, and valuable suggestions during the preparation of the paper. The author particularly thanks Jean-Michel Coron for the explanation of the dynamic feedback controls, Emmanuel Tr\'elat for extremely helpful information and comments on the linear quadratic optimal control theory and the theory of the Riccati equation,  and Marius Tusnack for the valuable information and helpful comments on the semigroup theory and its connection with the control theory, and the theory of the Riccati equation.

\providecommand{\bysame}{\leavevmode\hbox to3em{\hrulefill}\thinspace}
\providecommand{\MR}{\relax\ifhmode\unskip\space\fi MR }
% \MRhref is called by the amsart/book/proc definition of \MR.
\providecommand{\MRhref}[2]{%
  \href{http://www.ams.org/mathscinet-getitem?mr=#1}{#2}
}
\providecommand{\href}[2]{#2}

\end{document}